\documentclass[11pt]{paper}

\usepackage{amssymb,amsfonts,amsthm,amsmath}

\usepackage{latexsym,multicol,graphicx}

\usepackage{graphicx, verbatim}

\usepackage{xcolor}
\usepackage{verbatim}

\usepackage{pict2e}

\usepackage[a4paper, hmargin=2.6cm, vmargin={3.5cm, 3.5cm}]{geometry}

\usepackage{fancyhdr}
\usepackage[normalem]{ulem}

\usepackage[title,titletoc,toc]{appendix}

\newtheorem{theorem}{Theorem}
\newtheorem{lemma}{Lemma}

\renewcommand{\phi}{\varphi}

\newcommand{\e}{\varepsilon}

\newcommand{\al}{\alpha}

\newcommand{\eps}{\varepsilon}

\newcommand{\la}{\lambda}

\newcommand{\cE}{\mathcal E}

\newcommand{\bE}{\mathbb E}

\def\done{{1\hskip-2.5pt{\rm l}}}

\def\done{{1\hskip-2.5pt{\rm l}}}
\global\long\def\dd{\,{\rm d}}

\renewcommand{\le}{\leqslant}

\newcommand{\bR}{\mathbb R}
\newcommand{\bC}{\mathbb C}

\newcommand{\bD}{\mathbb D}
\newcommand{\bT}{\mathbb T}
\newcommand{\bN}{\mathbb N}

\newcommand{\bP}{\mathbb P}

\begin{document}

\title{Hole probability for zeroes of Gaussian Taylor series\\ with finite radii of convergence}

\author{Jeremiah Buckley
\thanks{Supported by ISF Grants~~1048/11 and~166/11, by ERC Grant~335141 and by the Raymond and Beverly Sackler Post-Doctoral
Scholarship 2013--14.}
\\
Department of Mathematics, King's College London, \\
Strand, London, WC2R 2LS, UK
\and
Alon Nishry \\
Department of Mathematics, University of Michigan, \\
530 Church Street, Ann Arbor, MI  48109, USA
\and
Ron Peled
\thanks{Supported by ISF Grants~1048/11 and~861/15 and by IRG Grant SPTRF.}\\
School of Mathematical Sciences, Tel Aviv University\\
Tel Aviv 69978, Israel
\and
Mikhail Sodin
\thanks{Supported by ISF Grants~166/11 and~382/15.}
\\
School of Mathematical Sciences, Tel Aviv University\\
Tel Aviv 69978, Israel
}

\date{\today}

\maketitle

\begin{abstract}
We study a family of random Taylor series
\begin{equation*}
  F(z) = \sum_{n\ge 0} \zeta_n a_n z^n
\end{equation*}
with radius of convergence almost surely $1$ and independent, identically
distributed complex Gaussian coefficients $(\zeta_n)$; these Taylor series are distinguished
by the invariance of their zero sets with respect to isometries of the
unit disk. We find reasonably tight upper and lower bounds on the
probability that $F$ does not vanish in the disk $\{|z|\le r\}$ as
$r\uparrow 1$. Our bounds take different forms according to whether
the non-random coefficients $(a_n)$ grow, decay or remain of the
same order. The results apply more generally to a class of Gaussian
Taylor series whose coefficients $(a_n)$ display power-law behavior.
\end{abstract}

\section{Introduction}\label{sect:intro}

Random analytic functions are a classical topic of study, attracting the attention of analysts and probabilists~\cite{Kahane}.
One of the most natural instances of such functions is provided by Gaussian analytic functions (GAFs, for short).
The zero sets of Gaussian analytic functions  are point processes exhibiting many interesting features~\cite{HKPV, NS1}.
These features depend on the geometry of the domain of analyticity of the function and, sometimes, pique the interest of
mathematical physicists, see, for instance,~\cite{BBL, FH, Hannay, NV}.

In this work we consider Gaussian analytic functions in the unit disk represented by a Gaussian Taylor series
\begin{equation}\label{eq1}
F(z) = \sum_{n\ge 0} \zeta_n a_n z^n\,,
\end{equation}
where $(\zeta_n)$ is a sequence of independent, identically
distributed standard complex normal random variables (i.e., the probability density of each $\zeta_n$ is
$\frac{1}{\pi}e^{-|z|^2}$ with respect to Lebesgue measure on the complex plane), and $(a_n)$ is a non-random sequence of
non-negative numbers satisfying $\limsup \sqrt[n]{a_n} = 1$. Properties of zero sets of Gaussian Taylor series with infinite radius of convergence have been studied quite intensively in
recent years, see, e.g.,~\cite{HKPV, NS, NishryT} and the references therein. When the radius of convergence is finite, the hyperbolic geometry often
leads to some peculiarities and complications. In particular, this is the case in the study of the hole probability, that is, the probability of the event
\begin{equation*}
\text{Hole}(r) = \bigl\{ F\ne 0 \ \text{on}\ r\bar\bD \bigr\},
\end{equation*}
an important characteristic both from the point of view of analytic function theory and the theory of point processes.
For \emph{arbitrary} Gaussian Taylor series with radius of convergence one,
understanding the asymptotic behaviour of the hole probability as $r\uparrow 1$ seems difficult. Here, we focus on a natural family $F_L$ of
hyperbolic Gaussian analytic functions (hyperbolic GAFs, for short),
\begin{equation}\label{eq9}
F_L(z) = \sum_{n\ge 0} \zeta_n \sqrt{\frac{L(L+1)\,\ldots\, (L+n-1)}{n!}}\, z^n\,, \qquad
0<L<\infty\,,
\end{equation}
distinguished by the invariance of the distribution of their zero set
under M\"obius transformations of the unit
disk~\cite[Ch.~2]{HKPV}, \cite{ST1}. Note that the parameter $L$ used to parameterize this family
equals the mean number of zeros of $F_L$ per unit hyperbolic area.

The main result of our work provides reasonably tight asymptotic bounds on the logarithm of the hole probability as $r$ increases to $1$.
These bounds show a transition in the asymptotics of the hole probability according to whether $0<L<1$, $L=1$ or $L>1$.
Curiously, a transition taking place at a different point, $L=\tfrac 12$,
was observed previously in the study of the asymptotics of the variance of the number of zeroes of $F_L$ in disks of increasing radii~\cite{Jerry}.

\subsection{The main result}\label{subsect:main-result}
\begin{theorem}\label{thm:main} Suppose that $F_L$ is a hyperbolic GAF,
and that $r\uparrow 1$. Then,

\smallskip\par\noindent {\rm (i)}
for $0<L<1$,
\[
\frac{1-L-o(1)}{2^{L+1}}\, \frac1{(1-r)^L}\, \log\frac1{1-r} \le -
\log\bP[{\rm Hole}(r)] \le  \frac{1-L+o(1)}{2^{L}}\,
\frac1{(1-r)^L}\, \log\frac1{1-r}\,;
\]

\medskip\par\noindent {\rm (ii)}
for $L=1$,
\[
-\log\bP[{\rm Hole}(r)] = \frac{\pi^2+o(1)}{12}\, \frac1{1-r}\,;
\]

\smallskip\par\noindent {\rm (iii)}
for $L>1$,
\[
-\log\bP[{\rm Hole}(r)] = \frac{(L-1)^2+o(1)}4\, \frac1{1-r}\, \log^2 \frac1{1-r}\,.
\]
\end{theorem}

\medskip\noindent
The case $L=1$ in this theorem is due to Peres and Vir\'ag. We'll briefly discuss their work in Section~\ref{subsubsect:Peres-Virag}. The rest is apparently new.

\subsection{Previous work}\label{subsect:previous-wors}

\subsubsection{Gaussian Taylor series with an infinite radius of
convergence}\label{sec:infinite_radius} For an {\em arbitrary}
Gaussian Taylor series with an \emph{infinite} radius of
convergence, the logarithmic asymptotics of the hole probability was
obtained in~\cite{Nishry}. The main result~\cite[Theorem 1]{Nishry} says
that when $r\to\infty$ outside an exceptional set of finite
logarithmic length,
\begin{equation}\label{eq4}
-\log\bP[{\rm Hole}(r)] = S_F(r) + o(S_F(r))\,,
\end{equation}
where
\begin{equation}\label{eq3}
S_F(r) = \sum_{n\ge 0} \log_+ (a_n^2 r^{2n})\,.
\end{equation}
In this generality, the appearance of an exceptional set of values
of $r$ is unavoidable due to possible irregularities in the
behaviour of the coefficients $(a_n)$ (see~\cite[Section 17]{NishryT}).

\medskip
For a Gaussian Taylor series with a finite radius of convergence the asymptotic rate of
decay of the hole probability has been described only in several rather special cases.

\subsubsection{The determinantal case: $L=1$}\label{subsubsect:Peres-Virag}
Peres and Vir\'ag~\cite{PV} (see also \cite[Section~5.1]{HKPV})
discovered that the zero set of the Gaussian Taylor
series
\begin{equation}\label{eq5}
F(z) = \sum_{n\ge 0} \zeta_n z^n
\end{equation}
(that corresponds to $L=1$ in \eqref{eq9}) is a determinantal point
process~\cite[Theorem 1]{PV} and, therefore, many of its characteristics can be
explicitly computed. In particular, they found that~\cite[Corollary 3 (i)]{PV}
\begin{equation}\label{eq6}
-\log\bP[{\rm Hole}(r)] = \frac{\pi^2+o(1)}{12}\, \frac1{1-r}\,\qquad\text{as }r\to 1.
\end{equation}
For $L\ne 1$, the zero set of $F_L$ is not a determinantal point process \cite[p.~83]{HKPV}, requiring the use of other techniques.

\subsubsection{Fast growing coefficients}
Skaskiv and Kuryliak~\cite{KS} showed that the technique developed in~\cite{Nishry}
can be applied to Gaussian Taylor series on the disk, that grow very fast near the boundary.
Put
\[
\sigma_F(r)^2 = \bE[|F(re^{{\rm i}\theta})|^2] = \sum_{n\ge 0} a_n^2 r^{2n}\,.
\]
They proved~\cite[Theorem 4]{KS} that if
\[
\lim_{r\to 1} (1-r) \log\log\sigma_F(r) = +\infty\,,
\]
and the sequence $(a_n)$ is logarithmically concave, then the same logarithmic
asymptotic~\eqref{eq4} holds when $r\to 1$ outside a small exceptional subset
of $[\tfrac12, 1)$. Note that in our case, $\sigma_{F_L}(r) = (1-r^2)^{-L/2}$.

\subsubsection{The case when $F$ is bounded on $\bD$}
At the opposite extreme, we may consider the case when, almost surely, the random series $F$ is
bounded on $\bD$. Then $\bP[\text{Hole}(r)]$ has a positive limit as $r\to 1$.
Indeed, put $F=F(0)+G$. If $F$ is bounded on $\bD$, then $G$ is bounded on $\bD$ as well.
Take $M$ so that $\bP[\sup_{\bD}|G| \le M]\ge \frac12$. Then
\[
\bigl\{ F\ne 0 \ \text{on}\ \bD \bigr\} \supset \bigl\{ \sup_{\bD}|G| \le M, \ |F(0)|>M \bigr\}\,.
\]
Since $F(0)$ and $G$ are independent, we get
\[
\bP[\text{Hole}(r)] \ge \bP[\sup_{\bD}|G| \le M] \cdot \bP[|F(0)|>M] \ge \frac12 e^{-M^2}\,.
\]

In view of this observation,
we recall a classical result that goes back to Paley and Zygmund and gives a sufficient condition
for continuity (and hence, boundedness) of $F$ on $\bar\bD$. Introduce the sequence
\[
s_j = \Bigl(\, \sum_{2^j \le n < 2^{j+1}} a_n^2 \Bigr)^{1/2}\,.
\]
If the sequence $s_j$ decreases and $\sum_j s_j < \infty$, then, almost surely,
$F$ is a continuous function on $\bar\bD$ \cite[Section~7.1]{Kahane}. On the other hand, under a mild
additional regularity condition, divergence of the series $\sum_j s_j = +\infty$ guarantees
that, almost surely, $F$ is unbounded in $\bD$ \cite[Section~8.4]{Kahane}.

\subsection{Several comments on Theorem~\ref{thm:main}}\label{subsect:brief-discussion}

\subsubsection{}
The proof of Theorem~\ref{thm:main} combines the tools introduced
in~\cite{ST3, Nishry} with several new ingredients. Unfortunately,
in the case $0< L <1$, our techniques are insufficient for finding
the main term in the logarithmic asymptotics of
$\bP[\text{Hole}(r)]$. Also we cannot completely recover the
aforementioned result of Peres and Vir\'ag. On the other hand, our
arguments do not use the hyperbolic invariance of the zero
distribution of $F_L$. We make use of the fact that
\begin{equation}\label{eq:a_n_asymptotics}
a_n^2 = \frac{\Gamma(n+L)}{\Gamma(L)
\Gamma(n+1)}=(1+o(1))\frac{n^{L-1}}{\Gamma(L)}\,,\quad n\to\infty,
\end{equation}
where $\Gamma$ is Euler's Gamma function, and our techniques apply
more generally to a class of Gaussian Taylor series whose
coefficients $(a_n)$ display power-law behavior. We will return to
this in the concluding Section~\ref{Section:Last}.

\subsubsection{}
In the case $L>1$, using \eqref{eq:a_n_asymptotics}, it is easy to
see that
\[
S_F(r) = \sum_{n\ge 0} \log_+ (a_n^2 r^{2n}) = \frac{(L-1)^2+o(1)}4\, \frac1{1-r}\, \log^2\frac1{1-r}\,,
\qquad r\to 1\,.
\]
That is, in this case the hyperbolic geometry becomes less relevant
and the main term in the logarithmic asymptotic of the hole
probability is governed by the same function as in the planar case,
discussed in Section~\ref{sec:infinite_radius}.

\subsubsection{}

For $0<L<1$, the gap between the upper and lower bounds in Theorem~\ref{thm:main} remains unsettled.
It would be also interesting to accurately explore the behaviour of the logarithm of the hole probability
near the transition points $L=0$ and $L=1$; for instance, to consider the cases $a_n = n^{-1/2}\ell (n)$ and
$a_n = \ell(n)$, where $\ell$ is a slowly varying function.

Of course, the ultimate goal would be to treat the case of {\em
arbitrary} Gaussian Taylor series with finite radii of convergence.

\subsubsection*{Acknowledgements}
The authors thank Alexander Borichev, Fedor Nazarov, and the referee for several useful suggestions.

\subsubsection*{Notation}

\begin{itemize}
\item Gaussian analytic functions will be called GAFs. The Gaussian analytic functions $F_L$ defined in~\eqref{eq9} will be called
hyperbolic GAFs.
\item We suppress all of the dependence on $L$ unless it is absolutely necessary. In particular,
from here on,
\begin{itemize}
\item[$\bullet$] $F=F_L$,
\item[$\bullet$] $c$ and $C$ denote positive constants that might \underline{only} depend on $L$. The values of
these constants are irrelevant for our purposes and may vary from
line to line. By $A$, $\al$, $\al_i$, etc. we denote positive
constants whose values we keep fixed throughout the proof in which
they appear.
\end{itemize}
\item
The notation $X\simeq Y$ means that $cX \le Y \le CY$.
\item We set
\begin{equation*}
  \delta = 1-r\quad\text{and}\quad r_0=1-\kappa\delta\quad\text{with $1<\kappa\le 2$}.
\end{equation*}
Everywhere, except Section~\ref{Sect:UB_L>1}, we set $\kappa=2$,
that is, $r_0=1-2\delta$. The value of $r$ is assumed to be
sufficiently close to $1$. Correspondingly, the value of $\delta$ is
assumed to be sufficiently small.
\item The variance $\sigma_F^2 = \sigma_F(r)^2$ is defined by
\[ \sigma_F^2 = \bE\bigl[ |F(re^{{\rm i}\theta})|^2\bigr] = (1-r^2)^{-L}
= (1+o(1)) (2\delta)^{-L}, \quad {\rm  as}\  \delta\to 0. \]
Usually, we suppress the dependence on $r$ and write
$\sigma_F$ instead of $\sigma_F(r)$. Notice that $\log\frac{1}{\delta} \simeq \log \sigma_F$.
\item An event $E$ depending on $r$ will be called \emph{negligible} if $-\log\bP[{\rm Hole}(r)] = o(1) \bigl( -\log\bP[E]] \bigr)$ as $r\to 1$.
Notice that this may depend on the value of $L$.
\item If $f$ takes real values, then we define $f_+ = \max\{0,f\}$ and $f_- = \max\{0,-f\}$.
\item $e(t)=e^{2\pi{\rm i}t}$.
\item $[x]$ denotes the integer part of $x$.
\item $n \equiv k\, (N)$ means that $n \equiv k$ modulo $N$.
\item $\bD$ denotes the open unit disk, $\bT$ denotes the unit circle.
\item The planar Lebesgue measure is denoted by $m$, and the (normalized) Lebesgue measure on $\bT$
is denoted by $\mu$.
\end{itemize}

\section{Idea of the proof}
We give a brief description of the proof of Theorem~\ref{thm:main}
in the cases $0<L<1$ and $L>1$. In the case $L=1$ our arguments
suffice to estimate the logarithm of the hole probability up to a
constant, as discussed in Section~\ref{Section:Last}, and we briefly
sketch the argument for this case as well. Our proof for the upper
bound on the hole probability in the case $0<L<1$ is more involved
than the other cases.
\subsection{Upper bounds on the hole probability when $L>1$ and $L=1$}
Our starting point for proving upper bounds is the mean-value
property. On the hole event,
\begin{equation*}
  \int_\bT \log|F(tr)| \dd \mu(t) = \log|F(0)|.
\end{equation*}
Off an event of negligible probability, the integral may be
discretized, yielding the inequality
\begin{equation}\label{eq:discrete_mean_value}
  \sum_{j=1}^N \log| F(\tau \omega^j r_0)| \le N \log |F(0)| + 1
\end{equation}
in the slightly smaller radius $r_0<r$, for a random $\tau \in \bT$, taken
out of a small set of possibilities, suitable $N$ and $\omega =
e(1/N)$ (see Lemma~\ref{lemma8} and Lemma~\ref{lemma9}). Thus it
suffices to bound from above, for each fixed $\tau\in\bT$, the
probability that \eqref{eq:discrete_mean_value} holds. We may
further simplify by fixing a threshold $T>0$, noting that
$\bP\left[|F(0)|\ge T\right] = \exp(-T^2)$ and writing
\begin{equation}\label{eq:after_threshold_value}
  \bP\Bigl[\sum_{j=1}^N \log| F(\tau \omega^j r_0)| \le N \log |F(0)| +
  1\Bigr]\le \bP\Bigl[\sum_{j=1}^N \log| F(\tau \omega^j r_0)| \le N \log T +
  1\Bigr] + e^{-T^2}.
\end{equation}
We focus on the first summand, setting $T$ sufficiently large so
that the second summand is negligible. Taking $0<\theta<2$ and
applying Chebyshev's inequality,
\begin{equation}\label{eq:exponential_Chebshev}
\begin{split}
\bP \Bigl[\sum_{j=1}^N \log| F(\tau \omega^j r_0)| \le N \log T + 1
\Bigr]
& \le \bP \Bigl[ \prod_{j=1}^N \bigl| F(\omega^j \tau r_0 ) \bigr|^{-\theta} \ge c T^{-\theta N} \Bigr] \\
& \le C T^{\theta N}\, \bE \Bigl[ \prod_{j=1}^N \bigl| F(\omega^j
\tau r_0 ) \bigr|^{-\theta}\Bigr]\,.
\end{split}
\end{equation}
It remains to estimate the expectation in the last expression. Our
bounds for it make use of the fact that the covariance matrix
$\Sigma$ of the Gaussian vector $(F(\tau \omega^j r_0))$, $1\le j\le
N$, has a circulant structure, allowing it to be explicitly
diagonalized. In particular, its eigenvalues are (see
Lemma~\ref{lemma14})
\begin{equation}\label{eq:eigenvalues}
\la_m = N\, \sum_{n \equiv m\, (N)} a_n^2 r^{2n}\,, \qquad m=0,
\ldots , N-1\,.
\end{equation}
This is used together with the following, somewhat rough, bound (see
Lemma~\ref{lemma18})
\begin{equation}\label{eq:determinant_bound}
\bE \Bigl[ \prod_{j=1}^N \bigl| F(\omega^j \tau r_0 )
\bigr|^{-\theta}\Bigr] \le \frac1{\det \Sigma}\, \bigg(
\Lambda^{\bigl(1-\tfrac12 \theta\bigl)} \cdot \Gamma \bigl(
1-\tfrac12\,\theta \bigr) \bigg)^N \,,
\end{equation}
where $\Lambda$ is the maximal eigenvalue of $\Sigma$ and $\Gamma$
is Euler's Gamma-function.

\subsubsection{The case $L>1$} We set the parameters to be $T =
\delta^{-\frac{1}{2}}\exp\left(-\sqrt{\log\frac{1}{\delta}}\right)$,
so that the factor $e^{-T^2}$ in \eqref{eq:after_threshold_value} is
indeed negligible, $N = \bigl[\,\frac{L-1}{2\delta}
\log\frac{1}{\delta}\,\bigr]$ and $\theta = 2 - \left(\log
\frac{1}{\delta}\right)^{-1}$. With this choice, the dominant term
in the combination of the bounds \eqref{eq:exponential_Chebshev} and
\eqref{eq:determinant_bound} is the factor $T^{\theta N} / \det
\Sigma$. Its logarithmic asymptotics are calculated using
\eqref{eq:eigenvalues} and yield the required upper bound.

\medskip
We mention that choosing $\theta$ close to its maximal value of $2$
corresponds, in some sense, to the fact that the event
\eqref{eq:discrete_mean_value} constitutes a very large deviation
for the random sum $\sum_{j=1}^N \log| F(\tau \omega^j r_0)|$.

\subsubsection{The case $L=1$} The same approach may be applied with
the parameters $T = b\,\delta^{-1/2}$, for a small parameter $b>0$,
$N = \bigl[\,\delta^{-1}\,\bigr]$ and $\theta = 1$. For variety,
Section~\ref{sec:upper_bound_L_1} presents a slightly different
alternative.

\medskip
In the same sense as before, the choice $\theta=1$ indicates that we
are now considering a large deviation event.

\subsection{Upper bound on the hole probability when $0<L<1$}
Our goal here is to show that the intersection of the hole event
with the event $\bigl\{|F(0)|\le A\sigma_F\sqrt{\log\tfrac{1}{(1-r)\sigma_F^2}}\,\bigr\}$ is
negligible when $A^2 < \frac{1}{2}$. The upper bound then follows
from the estimate
\[ \bP\Bigl[|F(0)|> A\sigma_F\sqrt{\log\tfrac{1}{(1-r)\sigma_F^2}}\,\Bigr] = \exp\left(-A^2\sigma_F^2\log\tfrac{1}{(1-r)\sigma_F^2}\right)\,. \]

The starting point is again the inequality
\eqref{eq:discrete_mean_value} in which we choose the parameter
\begin{equation*}
  N = \bigl[\delta^{-\alpha}\bigr],\quad L<\alpha<1,
\end{equation*}
with $\alpha$ eventually chosen close to $1$. However, a more
refined analysis is required here. First, we separate the constant
term from the function $F$, writing
\begin{equation*}
  F(z) = F(0) + G(z).
\end{equation*}
Second, to have better control of the Gaussian vector
\begin{equation}\label{eq:G_values}
(G(\tau \omega^j r_0)),\quad 1\le j\le N,
\end{equation}
we couple $G$ with two independent GAFs $G_1$ and $G_2$ so that $G =
G_1 + G_2$ almost surely, the vector $(G_1(\tau \omega^j r_0))$, $1\le j\le N$, is
composed of independent, identically distributed variables, and $G_2$
is a polynomial of degree $N$ with relatively small variance. In
essence, we are treating the variables in \eqref{eq:G_values} as
independent, identically distributed up to a small error captured by
$G_2$.

We proceed by conditioning on $F(0)$ and $G_2$, using the convenient
notation
\begin{align*}
\bE^{F(0), G_2} \bigl[\ . \ \bigr] &= \bE \bigl[\ . \ \big|\, F(0), G_2 \bigr], \\
\bP^{F(0), G_2} \bigl[\ . \ \bigr] &= \bP \bigl[\ . \ \big|\, F(0), G_2 \bigr].
\end{align*}
Applying Chebyshev's inequality we may use the above independence to
exchange expectation and product, writing, for $0<\theta<2$,
\begin{equation}\label{eq:Chebyshev_L_less_than_1}
\begin{split}
  \bP^{F(0), G_2}\Bigl[\sum_{j=1}^N \log| F(\tau \omega^j r_0)| \le N \log |F(0)| +
  1\Bigr] &\le C |F(0)|^{\theta N}\, \bE^{F(0), G_2} \Bigl[ \prod_{j=1}^N \bigl| F(\omega^j
\tau r_0 ) \bigr|^{-\theta}\Bigr]\\
&= C |F(0)|^{\theta N}\, \prod_{j=1}^N \bE^{F(0), G_2} \Bigl[ \bigl|
F(\omega^j \tau r_0 ) \bigr|^{-\theta}\Bigr].
\end{split}
\end{equation}
Thus we need to estimate expectations of the form
\begin{equation}\label{eq:negative_moments}
  \bE^{F(0), G_2} \Bigl[ \Bigl|\frac{F(\omega^j \tau r_0 )}{F(0)
} \Bigr|^{-\theta}\Bigr] = \bE^{F(0), G_2} \Bigl[ \Bigl|1 +
\frac{G_1(\omega^j \tau r_0 )+G_2(\omega^j \tau r_0 )}{F(0)}
\Bigr|^{-\theta}\Bigr],
\end{equation}
in which $F(0)$ and $G_2$ are given. Two bounds are used to this
end. Given a standard complex Gaussian random variable $\zeta$, real
$t>0$ and $0<\theta\le 1$ we have the simple estimate,
\begin{equation}\label{eq:simple_negative_moments}
\sup_{w\in\bC} \bE \Bigl[\, \Bigl| w+\frac{\zeta}{t} \Bigr|^{-\theta}  \, \Bigr]
\le t^\theta(1+C\theta)\,
\end{equation}
and, for $0\le \theta\le \frac{1}{2}$, the more refined
\begin{equation}\label{eq:precise_negative_moments}
\bE \Bigl[ \Bigl| 1+ \frac{\zeta}t \Bigr|^{-\theta} \Bigr]
\le 1 - c\theta\, \frac{e^{-t^2}}{1+t^2} + C\theta^2\,,
\end{equation}
see Lemma~\ref{lemma15} and Lemma~\ref{lemma17}. The most important feature of the second bound is that it is less than $1$ (though only slightly) when $\theta$ is very close to $0$,
satisfying $\theta\le c(1+t^2)^{-1}e^{-t^2}$.

Combining \eqref{eq:Chebyshev_L_less_than_1} and \eqref{eq:negative_moments} with the simple estimate \eqref{eq:simple_negative_moments} (with $\theta=1$) already suffices to prove that
the intersection of the hole
event with the event $\{|F(0)|\le a\sigma_F\}$ is negligible when
$a$ is a sufficiently small constant.
However, on the event
\[
\bigl\{ a\sigma_F\le |F(0)|\le A\sigma_F\sqrt{\log\tfrac{1}{(1-r)\sigma_F^2}}\, \bigr\}
\]
the error term $G_2$ becomes more
relevant and we consider two cases according to its magnitude.
Taking small constants $\e,\alpha_0>0$ and $\eta =
\delta^{\alpha_0}$ we let
\[
J = \Bigl\{1\le j \le N\colon \Bigl| 1+\frac{G_2(\omega^j
r_0)}{F(0)} \Bigr|\ge 1+2\e \Bigr\}\,.
\]
\subsubsection{The case $|J|>(1-2\eta)N$} Here, discarding (a priori,
before conditioning on $G_2$) a negligible event to handle the
rotation $\tau$, for many values of $j$, the terms in \eqref{eq:negative_moments}
have $\Bigl|1 + \frac{G_2(\omega^j \tau r_0 )}{F(0)}\Bigr|\ge 1+\e$.
This fact together with the bound
\eqref{eq:precise_negative_moments} (in simplified form, with
right-hand side $1 + C\theta^2$), taking $\theta$ tending to $0$ as
a small power of $\delta$, suffices to show that the probability in
\eqref{eq:Chebyshev_L_less_than_1} is negligible.

\subsubsection{The case $|J|\le(1-2\eta)N$} In this case we change
our starting point. By the mean-value inequality,
\begin{equation*}
  \int_\bT \log\Bigl|\frac{F(tr)}{F(0) + G_2(tr)}\Bigr| \dd \mu(t) \le \log\Bigl|\frac{F(0)}{F(0) + G_2(0)}\Bigr| = 0.
\end{equation*}
Off an event of negligible probability, the integral may again be
discretized, yielding
\begin{equation*}
  \sum_{j=1}^N \log\Bigl| \frac{F(\tau \omega^j r_0)}{F(0) + G_2(\tau \omega^j r_0)}\Bigr| \le 1
\end{equation*}
in the slightly smaller radius $r_0<r$, for a random $\tau$, taken
out of a small set of possibilities, see Lemma~\ref{lemma8}. As
before, for each fixed $\tau$, Chebyshev's inequality and
independence show that,
\begin{equation}\label{eq:Chebyshev_L_less_than_1_last_case}
\begin{split}
  \bP^{F(0), G_2}\Bigl[\sum_{j=1}^N \log\Bigl| \frac{F(\tau \omega^j r_0)}{F(0) + G_2(\tau \omega^j r_0)}\Bigr| \le 1\Bigr]
  &\le C\prod_{j=1}^N \bE^{F(0), G_2} \Bigl[ \Bigl| \frac{F(\tau
\omega^j r_0)}{F(0) + G_2(\tau \omega^j r_0)} \Bigr|^{-\theta}\Bigr]
\end{split}
\end{equation}
and we are left with the task of estimating terms of the form
\begin{equation}\label{eq:negative_moments2}
  \bE^{F(0), G_2} \Bigl[ \Bigl| \frac{F(\tau
\omega^j r_0)}{F(0) + G_2(\tau \omega^j r_0)} \Bigr|^{-\theta}\Bigr]
= \bE^{F(0), G_2} \Biggl[ \Biggl| 1+\frac{G_1(\tau \omega^j r_0)}{F(0)\left(1+
\frac{G_2(\tau \omega^j r_0)}{F(0)}\right)} \Biggr|^{-\theta}\Biggr].
\end{equation}
Again discarding (a priori) a negligible event to handle the rotation
$\tau$, for many values of $j$, the terms in \eqref{eq:negative_moments2} have $\Bigl|1
+ \frac{G_2(\omega^j \tau r_0 )}{F(0)}\Bigr|\le 1+\e$. These terms
are estimated by using \eqref{eq:precise_negative_moments} with $t\le
(1+4\e)A\,\sqrt{\log\tfrac{1}{(1-r)\sigma_F^2}}$. Correspondingly we set $\theta = c\eta
((1-r)\sigma_F^2)^{(1+10\e)A^2}$ and obtain that the probability in
\eqref{eq:Chebyshev_L_less_than_1_last_case} satisfies
\begin{equation*}
  \bP^{F(0), G_2}\Bigl[\sum_{j=1}^N \log\Bigl| \frac{F(\tau \omega^j r_0)}{F(0) + G_2(\tau \omega^j r_0)}\Bigr| \le
  1\Bigr]\le \exp(-c\eta^2((1-r)\sigma_F^2)^{2(1+10\e)A^2} N).
\end{equation*}
Recalling that $\sigma_F^2 = (1-r^2)^{-L}$, $\eta = \delta^{\alpha_0}$ and $N = [\delta^{-\alpha}]$, and choosing $\e$ and $\alpha_0$ close to
$0$ and $\alpha$ close to $1$ shows that this probability is negligible provided
that $A^2 < \frac{1}{2}$.

\medskip
In contrast to the cases $L>1$ and $L=1$, the fact that we take
$\theta$ tending to $0$ can be viewed as saying that we are now
considering a moderate deviation event.

\subsection{Lower bounds on the hole probability}
The proofs of our lower bounds on the hole probability are less involved than the proofs of the upper bounds and the reader is referred to the relevant sections for details.
We mention here that in all cases we rely on the same basic strategy: Fix a threshold $M>0$ and observe that
\begin{align}\label{eq:lower_bound}
\nonumber
  \bP\bigl[\text{Hole}(r)\bigr] &\ge \bP\Bigl[|F(0)|>M,\;\; \max_{r\bar \bD} |F - F(0)|\le M\Bigr] \\
  &= e^{-M^2}\cdot\bP\Bigl[\max_{r\bT} |F - F(0)|\le M\Bigr].
\end{align}

\subsubsection{The case $0<L<1$} Here we take
\[
M = \sqrt{1-L + 2\e} \cdot \sigma_F\, \sqrt{\log\tfrac{1}{1-r}}\,.
\]
To estimate the right-hand
side of \eqref{eq:lower_bound} we discretize the circle $r\bT$ into
$N = \bigr[(1-r)^{-(1+\e)}\bigr]$ equally-spaced points. We then use
Harg\'e's version of the Gaussian correlation inequality to
estimate, by bounding $F'$, the probability that the maximum attained on the circle $r\bT$ is not much bigger than the maximum
attained on these points and that the value $F$ attains at each of the points is
not too large.

\subsubsection{The case $L>1$} Here we take $M = \frac{1}{\sqrt{\delta}}\left(\log\frac{1}{\delta}\right)^{\alpha}$ for some $\frac{1}{2}<\alpha<1$. We also set
$N=\bigl[\frac{2 L}{\delta}\log\frac{1}{\delta}\bigr]$. We prove that the event on the right-hand side of \eqref{eq:lower_bound} becomes typical after
conditioning that the first $N$ coefficients in the Taylor series of $F$ are suitably small, and estimate the probability of the conditioning event.

\subsubsection{The case $L=1$} This case seems the most delicate of
the lower bounds. We take $M = B\sqrt{\frac{1}{1-r}}$ for a large
constant $B$. To estimate the probability on the right-hand side of
\eqref{eq:lower_bound} we write the Taylor series of $F$ as an
infinite sum of polynomials of degree $N =
\Bigl\lceil\frac{1}{1-r}\Bigr\rceil$ and use an argument based on
Bernstein's inequality and Harg\'e's version of the Gaussian
correlation inequality.

\section{Preliminaries}

Here, we collect several lemmas, which will be used in the proof
of Theorem~\ref{thm:main}.

\subsection{GAFs}

\begin{lemma}[\cite{HKPV}, Lemma~2.4.4]\label{lemma1}
Let $g$ be a GAF on $\bD$, and let $\sup_{\bD} \bE[|g|^2] \le \sigma^2$.
Then, for every $\la>0$,
\[
\bP \Big[\max_{\frac12\bar\bD} |g| > \la \sigma \Big] \le Ce^{-c\la^2}\,.
\]
\end{lemma}

\begin{lemma}\label{lemma2}
Let $f$ be a GAF on $\bD$, and $s \in (0,\delta)$. Put
\[
\sigma_f ^2(r) = \max_{r\bar\bD} \bE[|f|^2]\,.
\]
Then, for every $\la>0$,
\[
\bP \Bigl[ \max_{r\bar\bD} |f| > \la \sigma_f\bigl( r + s \bigr) \Bigr]
\le \frac{C}{s} \cdot e^{-c\la^2}\,.
\]
In particular, for every $\lambda>0$,
\[
\bP \Bigl[ \max_{r\bar\bD} |f| > \la \sigma_f\bigl( \tfrac12 (1+r) \bigr) \Bigr]
\le C \delta^{-1} \, e^{-c\la^2}\,.
\]
\end{lemma}

\noindent{\em Proof of Lemma~\ref{lemma2}}: Take an integer $N
\simeq \frac1s$ and consider the scaled functions
\[
g_j(w) = f(z_j+ s w), \qquad z_j = r e(j/N), \ j=1\,
\ldots\, N\,.
\]
Since $\max_{r\bar\bD} |f| = \max_{r\bT} |f(z)|$, the the first statement follows by applying
Lemma~\ref{lemma1} to each $g_j$ and using the
union bound. The second statement follows from the first by taking $s=\tfrac12 \delta= \tfrac12 (1-r)$.
\hfill $\Box$

\subsection{A priori bounds for hyperbolic GAFs}\label{subsect:a-priori-bounds}

\begin{lemma}\label{lemma3}
Suppose that $F$ is a hyperbolic GAF. Then, for $p>1$,
\[
\bP \Bigl[ \max_{(1-\frac12 \delta)\bar\bD} |F| \ge \sigma_F^p
\Bigr] \le C \delta^{-1} \,  e^{-c \sigma_F^{2p-2}}\,.
\]
\end{lemma}

\noindent{\em Proof of Lemma~\ref{lemma3}}:
Since $\sigma_F\big(1-\tfrac14 \delta\big) \le C \sigma_F$, this
follows from Lemma~\ref{lemma2}. \hfill $\Box$

\begin{lemma}\label{lemma4}
Suppose that $F$ is a hyperbolic GAF. Then
\[
\bP \Bigl[\, \max_{(1-2\delta)\bar\bD} |F| \le e^{-\log^2 \sigma_F} \,\Bigr]
\le C e^{-c\delta^{-1}\log^4\sigma_F}\,.
\]
\end{lemma}

\noindent{\em Proof of Lemma~\ref{lemma4}}:
Suppose that $ \max_{(1-2\delta)\bD} |F| \le e^{-\log^2 \sigma_F} $.
Then, by Cauchy's inequalities,
\[
|\zeta_n| a_n \le (1-2\delta)^{-n} e^{-\log^2 \sigma_F}\,,
\]
whence, for $n>0$, using the fact that $a_n \ge c n^{\frac12 (L-1)}$ (see~\eqref{eq:a_n_asymptotics}), we get
\[
|\zeta_n| \le C n^{\frac12 (1-L)} (1-2\delta)^{-n} e^{-\log^2 \sigma_F}
< C n^{\frac12} e^{c\delta n -\log^2 \sigma_F}\,.
\]
In the range $n\simeq \tfrac1{\delta}\log^2\sigma_F$, we get
\begin{equation}\label{eq:zeta-n}
|\zeta_n| \le C \bigl( \tfrac1{\delta} \log^2\sigma_F \bigr)^{\frac12} e^{-c\log^2 \sigma_F}
\le e^{-c\log^2\sigma_F}
\end{equation}
provided that $\delta$ is sufficiently small. The probability that
\eqref{eq:zeta-n} holds simultaneously for all such  $n$, does not exceed
\[
\Bigl( e^{-c\log^2\sigma_F} \Bigr)^{c\delta^{-1}\log^2\sigma_F} = e^{-c\delta^{-1}\log^4\sigma_F}\,,
\]
completing the proof. \hfill $\Box$

\medskip
Next, we define ``the good event'' $\Omega_{\tt g} = \Omega_{\tt g} (r)$ by
\begin{equation}\label{eq:good-event}
\Omega_{\tt g} = \bigl\{ \max_{(1-\frac12 \delta)\bar\bD} |F| \le \sigma_F^3  \bigr\}
\bigcap \bigl\{ \max_{(1-2\delta)\bar\bD} |F| \ge e^{-\log^2 \sigma_F}  \bigr\}
\end{equation}
and note that by Lemma~\ref{lemma3} and Lemma~\ref{lemma4} the event
$\Omega_{\tt g}^c$ is negligible.

\begin{lemma}\label{lemma5}
Suppose that $F$ is a hyperbolic GAF. If $\gamma > 1$ then
\[
{\rm Hole}(r) \bigcap \Omega_{\tt g}(r) \subset \bigl\{\,  \min_{(1-\gamma \delta)\bar\bD}\, |F|  \ge
\exp\bigl[- C (\gamma - 1)^{-1} \delta^{-3} \bigr]   \, \bigr\}\,.
\]
\end{lemma}

\noindent{\em Proof of Lemma~\ref{lemma5}}: Suppose that we are on the event $ {\rm Hole}(r) \bigcap \Omega_{\tt g}(r) $.
Since $F$ does not vanish in $r\bD$, the function $\log |F|$ is harmonic therein, and therefore,
\begin{align*}
\max_{(1-\gamma \delta)\bar\bD}\, \log\frac{1}{|F(z)|}
& = \max_{(1-\gamma \delta)\bar\bD}\, \int_\bT \frac{r^2 - |z|^2}{|rt - z|^2} \log\frac{1}{|F(rt)|} \dd \mu(t) \\ \\
&\le \max_{(1-\gamma\delta)\bar\bD}\, \Bigl( \frac{r +|z|}{r-|z|} \Bigr)  \int_\bT \left| \log|F(rt)| \right| \dd \mu(t) \\ \\
& \le \frac{2}{(\gamma - 1) \delta} \int_\bT \left| \log|F(rt)| \right| \dd \mu(t). \\
\end{align*}
Furthermore, on $\Omega_{\tt g}$, let $w\in (1-2\delta)\bT$ be a point where $|F(w)|\ge e^{-\log^2\sigma_F}$. Then
\begin{align*}
-\bigl( \log^2\sigma_F \bigr) \le \log|F(w)|
& = \int_\bT \frac{r^2 - |w|^2}{|rt-w|^2}\, \log |F(rt)|\, {\rm d}\mu(t) \\ \\
& = \int_\bT \frac{r^2 - |w|^2}{|rt-w|^2}\,
\bigl[ \log_+ |F(rt)| - \log_- |F(rt)|]\, {\rm d}\mu(t)\,,
\end{align*}
whence,
\begin{align*}
\int_{\bT} \log_-|F(rt)|\, {\rm d}\mu(t) & \le
\frac{r+|w|}{r-|w|}\, \int_{\bT} \frac{r^2-|w|^2}{|rt-w|^2}\, \log_- |F(rt)|\, {\rm d}\mu(t) \\ \\
& \le \frac{r+|w|}{r-|w|}\, \Bigl[ \int_{\bT} \frac{r^2-|w|^2}{|rt-w|^2}\, \log_+ |F(rt)|\, {\rm d}\mu(t) + \log^2\sigma_F \Bigr] \\ \\
& \le \frac{r+|w|}{r-|w|} \, \Bigl[ \max_{r\bT} \log |F| +  \log^2\sigma_F \Bigr] \\ \\
& \le \frac2{\delta} \bigl[\, 3 \log\sigma_F +  \log^2\sigma_F\, \bigl]
\le \frac{C}{\delta} \cdot  \log^2 \sigma_F\,.
\end{align*}
Then,
\[
\max_{(1-\gamma \delta)\bar\bD}\, \log\frac{1}{|F|}
\le \frac{2}{(\gamma - 1) \delta} \int_{\bT} \bigl| \log |F(rt)| \bigr|\, {\rm d}\mu(t)
\le \frac{C}{(\gamma - 1) \delta^2} \cdot \log^2 \sigma_F
< \frac{C}{(\gamma - 1) \delta^3}\,,
\]
proving the lemma. \hfill $\Box$

\subsection{Averaging $\log |F|$ over roots of unity}

We start with a polynomial version.

\begin{lemma}\label{lemma6}
Let $S$ be a polynomial of degree $n \ge 1$, let $k\ge 4$ be an integer, and
let $\omega=e(1/k)$.
Then there exists $\tau$ with $\tau^{k^2 n}=1$ so that
\[
\frac1{k}\, \sum_{j=1}^k \log |S(\tau\omega^j)| \ge \log |S(0)| - \frac{C}{k^2}\,.
\]
\end{lemma}

\noindent{\em Proof of Lemma~\ref{lemma6}}:
Assume that $S(0)=1$ (otherwise, replace $S$ by $S/S(0)$). Then
$ S(z) = \prod_{\ell=1}^n (1-s_\ell z) $ and
\[
\prod_{j=1}^k S(\omega^j z) = \prod_{\ell=1}^n (1-s_\ell^k z^k) = S_1(z^k)\,,
\]
where
$ S_1(w) \stackrel{\rm def}= \prod_{\ell=1}^n (1-s_\ell^k w) $
is a polynomial of degree $n$.

Let $M=\max_{\bT} |S_1|$. By the maximum principle,
$M \ge |S_1(0)|=1$. By Bernstein's inequality,
$\max_{\bT} |S_1'|\le nM$.
Now, let $t_0=e^{{\rm i}\phi_0}\in\bT$ be a point where $|S_1(t_0)|=M$. Then, for any $t=e^{{\rm i}\phi}\in\bT$ with
$|\phi-\phi_0|< \tfrac{\eps}n$ (with $0\le \e<1$), we have
\[
|S_1(t)| \ge M - \frac{\eps}n \cdot nM = (1-\eps)M \ge 1-\eps\,.
\]
The roots of unity of order $k n$ form a $\frac{\pi}{k n}$-net on $\bT$. Thus,
there exists $t$ such that $t^{k n}=1$ and $\log |S_1(t)| \ge - \tfrac{C}{k}$.
Taking $\tau$ so that $\tau^k=t$, we get
\[
\frac1{k}\, \sum_{j=1}^k \log |S(\tau\omega^j)| = \frac1{k}\, \log |S_1(\tau^k)| \ge - \frac{C}{k^2}\,,
\]
while $\tau^{k^2 n}=1$. \hfill $\Box$

\begin{lemma}\label{lemma7}
Let $f$ be an analytic function on $\bD$ such that $M\stackrel{\rm
def}=\sup_{\bD} |f|<+\infty$. Let $0<\rho <1$, and denote by $q$ the
Taylor polynomial of $f$ around $0$ of degree $N$. If
\begin{equation}\label{eq:N}
N \ge \frac{1}{1-\rho}\, \log\frac{M}{1-\rho}\,,
\end{equation}
then $\max_{\rho\bar\bD} |f-q|<1$.
\end{lemma}

\noindent{\em Proof of Lemma~\ref{lemma7}}:
Let
\[
f(z) = \sum_{n\ge 0} c_n z^n\,.
\]
Then, by Cauchy's inequalities, $|c_n|\le M$, whence
\begin{align*}
\max_{\rho\bar\bD} |f-q| \le \sum_{n\ge N+1} |c_n| \rho^n & \le M\,
\frac{\rho^{N+1}}{1-\rho}=\frac{M}{1-\rho}\,(1-(1-\rho))^{N+1}<\frac{M}{1-\rho}\,e^{-N(1-\rho)}\le
1.\hfill\qed
\end{align*}

\begin{lemma}\label{lemma8}
Let $0 < m < 1$, $0< \rho < 1$, $k\ge 4$ be an integer, and $\omega = e(1/k)$. Suppose that $F$ is an analytic function on $\bD$ with $m\le \inf_{\bD}|F|  \le\sup_{\bD} |F|  \le m^{-1}$
and that $P$ is a polynomial with $P(0) \ne 0$.
There exist a positive integer
\[
K \le k^2 \Bigl[ \deg P + \frac{C}{1-\rho}\, \log\frac{k}{m(1-\rho)} \Bigr]\,,
\]
and $\tau \in \bT$ satisfying $\tau^{K}=1$
such that
\[
\frac1k\, \sum_{j=1}^k \log\Bigl| \frac{F}{P}(\tau\omega^j \rho) \Bigr|
\le \log\Bigl| \frac{F}{P}(0) \Bigr| + \frac{C}{k^2}\,.
\]
\end{lemma}

\noindent{\em Proof of Lemma~\ref{lemma8}}:
Let $0 < \eps < m$ be a small parameter. Applying Lemma~\ref{lemma7} to the function $f=(\eps F)^{-1}$, we get
a polynomial $Q$ with
\[
\deg Q \le \frac{C}{1-\rho}\, \log\frac1{\eps m (1-\rho)}\,,
\]
such that
\[
\max_{\rho\bar\bD} \Bigl| \frac1{F} - Q \Bigr| < \eps\,.
\]
Then, assuming that $\eps<\frac12 m$, we get
\begin{equation}\label{eq:A}
\max_{\rho\bar\bD} \bigl| \log |Q| + \log |F| \bigr|
= \max_{\rho\bar\bD} \bigl| \log | 1 + F(Q-\tfrac1{F})|\, \bigr|
\le C\eps m^{-1}\,.
\end{equation}
Applying Lemma~\ref{lemma6} to the polynomial $S=P\cdot Q$
and taking into account~\eqref{eq:A}, we see that there exists $\tau$ so that $\tau^{k^2\deg S}=1$ and
\begin{align*}
\frac1k\, \sum_{j=1}^k \log\Bigl| \frac{F}{P}(\tau\omega^j \rho) \Bigr|
& \le - \frac1k\, \sum_{j=1}^k \log\Bigl| S(\tau\omega^j \rho) \Bigr| + C\eps m^{-1} \\ \\
& \le - \log|S(0)| + C\bigl[ k^{-2} + \eps m^{-1} \bigr] \\ \\
& \le \log\Bigl| \frac{F}{P}(0) \Bigr| + C\bigl[ k^{-2} + \eps m^{-1} \bigr]\,.
\end{align*}
It remains to let $\eps= \tfrac12 m k^{-2} $ and
$ K = k^2\deg S = k^2 (\deg P + \deg Q)$. \hfill $\Box$

\medskip
We will only need the full strength of this lemma in Section~\ref{subsubsect:smallJ}. Elsewhere, the following lemma will be sufficient.
\begin{lemma}\label{lemma9}
Let $F$ be a hyperbolic GAF. Let
\[
1 < \kappa \le  2, \quad r_0 = 1-\kappa \delta\,.
\]
Let $k\ge 4$ be an integer and $\omega=e(1/k)$. Then
\[
\bP\bigl[ {\rm Hole}(r) \bigr] \le \frac{ck^2\log k}{(\kappa - 1)^2 \delta^{4}} \cdot \sup_{\tau\in\bT}
\bP\Bigl[
\sum_{j=1}^k \log| F(\tau \omega^j r_0)| \le k \log |F(0)| + C
\Bigr] + \langle {\rm negligible\ terms} \rangle\,.
\]
\end{lemma}
\noindent{\em Proof of Lemma~\ref{lemma9}}: Outside the negligible
event $ {\rm Hole}(r)\setminus \Omega_{\tt g}$ we have the bound
\[
\max_{(1-\frac12 \delta)\bar\bD} |F| \le \sigma_F^3 \le \frac{C}{\delta^{3L/2}}.
\]
According to Lemma~\ref{lemma5} applied with $\gamma=1+\tfrac12 (\kappa-1)$, we have
\[
{\rm Hole}(r)\setminus \Omega_{\tt g} \subset \Bigl\{ \min_{(1-\gamma \delta)\bar\bD}\, |F| \ge \exp\bigl[ -C (\kappa - 1)^{-1} \delta^{-3} \bigr]\Bigr\}\,.
\]
Therefore, letting $m= \exp\bigl[-C (\kappa - 1)^{-1} \delta^{-3} \bigr]$, we get
\[
\bP\big[ {\rm Hole}(r) \big] \le \bP\Bigl[
m \le \min_{(1-\gamma \delta)\bar\bD}\, |F| \le \max_{(1-\gamma \delta)\bar\bD}\, |F| \le m^{-1}
\Bigr] + \langle {\rm negligible\ terms} \rangle\,.
\]
Applying Lemma~\ref{lemma8} to the function $F\bigl( (1-\gamma\delta) z \bigr)$ with $P \equiv 1$ and
\[
\rho = \frac{r_0} {(1-\gamma\delta)} = 1 - \frac{(\kappa-\gamma)\delta}{1-\gamma\delta}\,,
\]
we find that
\[
\Bigl\{ m \le \min_{(1-\gamma \delta)\bar\bD}\, |F| \le \max_{(1-\gamma \delta)\bar\bD}\, |F| \le m^{-1} \Bigr\}
\subset \bigcup_{\tau\colon \tau^K=1} \Bigl\{  \sum_{j=1}^k \log\Bigl| F(\tau\omega^j r_0) \Bigr|
\le k \log|F(0)| + \frac{C}{k}\Bigr\}
\]
with some
\[
K \le \frac{Ck^2}{(\kappa-\gamma)\delta}\, \log\frac{k}{m(\kappa-\gamma)\delta} \le \frac{Ck^2\log k}{(\kappa-1)^2\delta^4}\,,
\]
completing the proof. \hfill $\Box$

\medskip
Note that everywhere except for Section~\ref{Sect:UB_L>1} we use this lemma with $\kappa=2$.

\subsection{The covariance matrix}

We will constantly exploit the fact that the covariance matrix of
the random variables $ F(\omega^j r \bigr)$, $\omega = e(1/N)$,
$1\le j \le N$, has a simple structure. This is valid for general
Gaussian Taylor series, as the following lemma details.

\begin{lemma}\label{lemma14}
Let $F$ be any Gaussian Taylor series of the form \eqref{eq1} with
radius of convergence $R$ and let $z_j = re(j/N)$ for $r<R$ and
$j=0, \ldots , N-1$. Consider the covariance matrix
$\Sigma = \Sigma(r, N)$ of the random variables $F(z_0)$, \ldots , $F(z_{N-1})$,
that is,
\[
\Sigma_{jk} = \bE \Bigl[ F(z_j) \overline{F(z_k) }\Bigr]
= \sum_{n\ge 0} a_n^2 r^{2n} e((j-k)n/N)\,.
\]
Then, the eigenvalues of $\Sigma$ are
\[
\la_m = N\, \sum_{n \equiv m\, (N)} a_n^2 r^{2n}\,, \qquad m=0,
\ldots , N-1\,,
\]
where $n\equiv m\, (N)$ denotes that $n$ is equivalent to $m$ modulo
$N$.
\end{lemma}

\noindent{\em Proof of Lemma~\ref{lemma14}}:
Observe that
\[
F(z_j) = \sum_{m=0}^{N-1}\, \sum_{n \equiv m\, (N)}\, \zeta_n a_n r^n e(jn/N)
=\sum_{m=0}^{N-1}\, e(jm/N)\, \sum_{n\equiv m\, (N)}\, \zeta_n a_n r^n\,.
\]
Define the $N \times N$ matrix $U$ by
\[
U_{jm} = \frac1{\sqrt{N}}\, e(jm/N)\,, \qquad 0 \le j, m \le N-1\,,
\]
and the vector $\Upsilon$ by
\[
\Upsilon_m = \sqrt{N}\, \sum_{n\equiv m\, (N)}\, \zeta_n a_n r^n\,.
\]
Then $U$ is a unitary matrix (it is the discrete Fourier transform matrix) and the
components of $\Upsilon$ are independent complex Gaussian random variables with
\[
\bE \bigl[\, |\Upsilon_m|^2\, \bigr] = N \sum_{n\equiv m\, (N)}\, a_n^2 r^{2n}\,.
\]
Finally, note that
\[
\begin{pmatrix}
F(z_0) \\
F(z_1) \\
\vdots \\
F(z_{N-1})
\end{pmatrix}
= U \Upsilon\,.
\]
Hence, the covariance matrices of $ \bigl( F(z_0)$, \ldots , $F(z_{N-1}) \bigr)$ and of $\Upsilon$
have the same set of eigenvalues. \hfill $\Box$

\subsection{Negative moments of Gaussian random variables}

In the following lemmas $\zeta$ is a standard complex Gaussian random
variable. Recall that, for $\theta<2$,
\begin{equation}\label{eq24}
\bE\bigl[\, |\zeta|^{-\theta}\, \bigr] = \Gamma\bigl( 1-\tfrac12\, \theta \bigr)\,,
\end{equation}
where $\Gamma$ is Euler's Gamma-function.

\begin{lemma}\label{lemma15}
There exists a numerical constant $C>0$ such that, for every $t>0$
and $0<\theta\le 1$,
\[
\sup_{w\in\bC} \bE \Bigl[\, \Bigl| w+\frac{\zeta}{t} \Bigr|^{-\theta}  \, \Bigr]
\le t^\theta(1+C\theta)\,.
\]
\end{lemma}

\noindent{\em Proof of Lemma~\ref{lemma15}}:
Write
\[
\bE \Bigl[\, \Bigl| w+\frac{\zeta}{t} \Bigr|^{-\theta}  \, \Bigr]
= \frac1\pi\, \int_{\bC} \Bigl| w+\frac{z}{t} \Bigr|^{-\theta} e^{-|z|^2}\, {\rm d}m(z)\,,
\]
where $m$ is the planar Lebesgue measure,
and use the Hardy-Littlewood rearrangement inequality, noting that the symmetric decreasing
rearrangement of $ \bigl| w+\tfrac{z}{t} \bigr|^{-\theta} $ is $ \bigl| \frac{z}t\bigr|^{-\theta} $,
and that $e^{-|z|^2}$ is already symmetric and decreasing. \hfill $\Box$

\begin{lemma}\label{lemma16}
For each $\tau>0$ there exists a $C(\tau)>0$ such that for every
integer $n\ge 1$,
\[
\sup_{|t|\ge \tau} \bE \Bigl[\, \Bigl| \log\bigl| 1+ \frac{\zeta}t \bigr|\, \Bigr|^n \, \Bigr]
\le C(\tau) n!\,.
\]
\end{lemma}

\noindent{\em Proof of Lemma~\ref{lemma16}}:
By the symmetry of $\zeta$ we may assume that $t > 0$. Put $X=\bigl| 1+\tfrac{\zeta}t \bigr|$ and write
\[
\bE \bigl[ \bigl| \log X \bigr|^n \bigr] = \bE \bigl[ \bigl(
\log\tfrac1X \bigr)^n \done_{\{X\le\frac12\} }\bigr] + \bE \bigl[
\bigl( \log\tfrac1X \bigr)^n \done_{\{\frac12< X\le 1\} } \bigr] +
\bE \bigl[ \bigl( \log X \bigr)^n \done_{\{X > 1\}}  \bigr]\,.
\]
We have
\begin{align*}
\bE \bigl[ \bigl( \log\tfrac1X \bigr)^n \done_{\{X\le\frac12\}
}\bigr] & = \frac1{\pi}\, \int_{\{ | 1+\frac{z}t|\le\frac12\}}
\Bigl( \log \bigl| 1+\frac{z}t \bigr|^{-1} \Bigr)^n
e^{-|z|^2}\, {\rm d}m(z) \\ \\
& \le \frac{e^{-t^2/4}}{\pi}\,
\int_{\{| 1+\frac{z}t|\le\frac12\}} \Bigl( \log \bigl| 1+\frac{z}t \bigr|^{-1} \Bigr)^n\, {\rm d}m(z) \\ \\
& =  \frac{e^{-t^2/4}t^2}{\pi}\,
\int_{\{|w|\le\frac12\}} \Bigl( \log\frac1{|w|}\Bigr)^n\, {\rm d}m(w) \\ \\
& \le C\, \int_0^{1/2} \Bigl( \log\frac1{s} \Bigr)^n s\,{\rm d}s \\ \\
& = C \int_{\log 2}^\infty x^n e^{-2x}\, {\rm d}x
\le C n!\,.
\end{align*}
In addition,
\[
\bE \bigl[ \bigl( \log\tfrac1X \bigr)^n \done_{\{\frac12< X\le 1\} } \bigr]
\le \bigl( \log 2 \bigr)^n < 1\,.
\]
Finally, for $t\ge \tau$,
\[
\bE \bigl[ \bigl( \log X \bigr)^n \done_{\{X > 1\}}  \bigr]
\le \bE \bigl[ (X-1)_+^n \bigr] \le \bE \bigl[ |X-1|^n \bigr]
\le \bE \Bigl[ \frac{|\zeta|^n}{t^n} \Bigr]
\stackrel{\eqref{eq24}}= \frac{\Gamma (\frac12 n + 1)}{t^n} \le \tau^{-n}n!\,,
\]
completing the proof. \hfill $\Box$

\begin{lemma}\label{lemma16.5}
For $t>0$,
\[
\bE \Bigl[ \log \Bigl| 1+ \frac{\zeta}t \Bigr| \Bigr] > \frac{e^{-t^2}}{2(t^2+1)}\,.
\]
\end{lemma}

\noindent{\em Proof of Lemma~\ref{lemma16.5}}: We have
\begin{align*}
\bE \Bigl[ \log \Bigl| 1+ \frac{\zeta}t \Bigr| \Bigr]
& = \frac1{\pi}\, \int_{\bC} \log \Bigl| 1+ \frac{z}t \Bigr|\, e^{-|z|^2}\, {\rm d}m(z) \\ \\
& = 2 \int_0^\infty e^{-r^2} r\,
\int_{-\pi}^\pi \log\Bigl| 1+\frac{re^{{\rm i}\theta}}t \Bigr|\, \frac{{\rm d}\theta}{2\pi} \, {\rm d}r\\ \\
& = 2 \int_0^\infty \log_+ \Bigl( \frac{r}{t} \Bigr)\,  e^{-r^2} r\, {\rm d}r \\ \\
& = \frac12\, \int_{t^2}^\infty \frac{e^{-u}}{u}\, {\rm d}u = I\,.
\end{align*}
Integrating by parts once again, we see that
\[
I = \frac{e^{-t^2}}{2t^2} - \int_{t^2}^\infty \frac{e^{-u}}{2u^2}\, {\rm d}u >  \frac{e^{-t^2}}{2t^2}
- \frac1{t^2}\, I\,,
\]
whence,
\[
I > \frac{e^{-t^2}}{2(t^2+1)}\,,
\]
completing the proof. \hfill $\Box$

\begin{lemma}\label{lemma17}
There exist numerical constants $c, C > 0$ such that, for every
$t>0$ and $0\le\theta\le \tfrac12$,
\[
\bE \Bigl[ \Bigl| 1+ \frac{\zeta}t \Bigr|^{-\theta} \Bigr]
\le 1 - c\theta\, \frac{e^{-t^2}}{1+t^2} + C\theta^2\,.
\]
\end{lemma}

\noindent{\em Proof of Lemma~\ref{lemma17}}:
Lemma~\ref{lemma15} yields that there exist $c, \tau>0$ such that,
for every $0<t\le \tau$ and $0\le\theta\le\tfrac12$, we have
\[
\bE \Bigl[ \Bigl| 1+ \frac{\zeta}t \Bigr|^{-\theta} \Bigr]
\le 1 - c\theta\,.
\]
Thus, we need to consider the case $t\ge \tau$. Write
\[
\Bigl| 1+ \frac{\zeta}t \Bigr|^{-\theta} = \exp\left(-\theta\log\Bigl|
1+ \frac{\zeta}t \Bigr|\right) = 1 + \sum_{n\ge 1} \frac{\bigl(
-\theta \log|1+\frac{\zeta}t |\bigr)^n}{n!}\,.
\]
Using Lemma~\ref{lemma16}, we get
\begin{multline*}
\bE \Bigl[ \Bigl| 1+ \frac{\zeta}t \Bigr|^{-\theta} \Bigr] \le 1 -
\theta\, \bE \Bigl[ \log \Bigl| 1+ \frac{\zeta}t \Bigr| \Bigr]
+ \sum_{n\ge 2} \frac{\theta^n}{n!}\, \bE \Bigl[ \Bigl| \log \bigl| 1+ \frac{\zeta}t \bigr|\Bigr|^n \Bigr] \\
\stackrel{0\le\theta\le \frac12}\le 1 - \theta\, \bE \Bigl[ \log
\Bigl| 1+ \frac{\zeta}t \Bigr| \Bigr] +
C(\tau)\,\frac{\theta^2}{1-\theta} \le 1 - \theta\, \bE \Bigl[ \log
\Bigl| 1+ \frac{\zeta}t \Bigr| \Bigr] + 2C(\tau)\theta^2\,.
\end{multline*}
Applying Lemma~\ref{lemma16.5}, we get the result. \hfill $\Box$

\begin{lemma}\label{lemma18}
Suppose that $(\eta_j)_{1\le j \le N}$ are complex Gaussian
random variables with covariance matrix $\Sigma$. Then,
for $0\le\theta < 2$,
\[
\bE \Bigg[ \prod_{j=1}^N \frac1{|\eta_j|^\theta}\,\Bigg]
\le \frac1{\det \Sigma}\, \bigg( \Lambda^{\bigl(1-\tfrac12 \theta\bigl)} \cdot \Gamma \bigl( 1-\tfrac12\,\theta \bigr) \bigg)^N \,,
\]
where $\Lambda$ is the maximal eigenvalue of $\Sigma$.
\end{lemma}

\noindent{\em Proof of Lemma~\ref{lemma18}}: We have
\begin{align*}
\bE \Bigl[ \prod_{j=1}^N \frac1{|\eta_j|^\theta}\,\Bigr] & =
\frac1{\pi^N \det\Sigma}\, \int_{\bC^N} \prod_{j=1}^N
\frac1{|Z_j|^\theta}\,
e^{-\langle\Sigma^{-1}Z, Z\rangle}\, {\rm d}m(Z_1) \ldots {\rm d}m(Z_N)  \\ \\
& \le \frac1{\pi^N \det\Sigma}\, \int_{\bC^N} \prod_{j=1}^N
\frac1{|Z_j|^\theta}\,
e^{-\Lambda^{-1} |Z|^2}\, {\rm d}m(Z_1) \ldots {\rm d}m(Z_N) \\ \\
& = \frac1{\det\Sigma}\, \Bigl( \frac1{\pi}\int_{\bC} |z|^{-\theta}\, e^{-\Lambda^{-1}|z|^2}\, {\rm d}m(z) \Bigr)^N \\ \\
& \stackrel{\eqref{eq24}}= \frac1{\det \Sigma}\, \bigg( \Lambda^{\bigl(1-\tfrac12 \theta\bigl)} \cdot \Gamma \bigl( 1-\tfrac12\,\theta \bigr) \bigg)^N\,,
\end{align*}
proving the lemma. \hfill $\Box$

\section{Upper bound on the hole probability for $0< L < 1$}
Now we are ready to prove the lower bound part of Theorem~\ref{thm:main}, in the case $0 < L < 1$.
Throughout the proof, we use the parameters
\[
L<\alpha<1, \quad r_0 = 1-2\delta, \quad N=\bigl[ \delta^{-\alpha} \bigr], \quad \omega = e(1/N),
\]
where
\[
\frac12 \le r < 1, \quad \delta = 1 - r.
\]
In many instances, we assume that $r$ is sufficiently close to $1$, that is, that $\delta$ is
sufficiently small.

\medskip
It will be convenient to separate the constant term from the
function $F$, letting
\begin{equation*}
F=F(0)+G.
\end{equation*}

\subsection{Splitting the function $G$}\label{sec:splitting_G}

We define two independent GAFs $G_1$ and $G_2$ so that $G=G_1+G_2$ and
\begin{itemize}
\item $G_1(\omega^j r_0)$, $1\le j \le N$, are independent, identically distributed Gaussian random variables with variance close to $\sigma_F^2(r_0)$;
\item $G_2$ is a polynomial of degree $N-1$ and, for $|z|=r_0$, the variance of $G_2(z)$ is much smaller than
$\sigma_F^2(r_0)$.
\end{itemize}
Let $(\zeta_n')_{n\ge 1}$ and $(\zeta_n'')_{1\le n \le N-1}$ be two independent sequences of
independent standard complex Gaussian random variables, and let
\begin{align*}
G_1(z) & = \sum_{n=1}^\infty \zeta_n' b_n z^n\,, \\
G_2(z) & = \sum_{n=1}^{N-1} \zeta_n'' d_n z^n\,,
\end{align*}
where the non-negative coefficients $b_n$ are defined by
\[
\begin{cases}
b_{n}^2 r_{0}^{2n}=\displaystyle{\sum_{k\ge1}}\Bigl[a_{kN}^2 r_0^{2kN}-a_{kN+n}^{2}r_0^{2(kN+n)}\Bigr],\quad & 1\le n\le N-1,\\
b_{n}=a_{n}, & n\ge N.
\end{cases}
\]
Since the sequence $(a_n)$ does not increase, the expression in the brackets is positive.
For the same reason, for $1\le n \le N-1$, we have
\[
a_n^2 r_0^{2n} + \sum_{k\ge 1} a_{kN+n}^2 r_0^{2(kN+n)} = \sum_{k\ge
0} a_{kN+n}^2 r_0^{2(kN+n)} \ge  \sum_{k\ge 1} a_{kN}^2 r_0^{2kN}\,,
\]
whence, for these values of $n$ we have $b_n\le a_n$. The
coefficients $d_n \ge 0$ are defined by $ a_n^2 = b_n^2 + d_n^2$.
This definition implies that the random Gaussian functions $ G $ and
$ G_1 + G_2 $ have the same distribution, and we couple $G$, $G_1$
and $G_2$ (that is, we couple the sequences $(\zeta_n)$,
$(\zeta_n')$ and $(\zeta_n'')$) so that $G =  G_1+G_2$  almost
surely.

\begin{lemma}\label{lemma18a}
For any $\tau\in\bT$, the random variables $\bigl( G_1(\omega^j \tau
r_0\bigr) \bigl)$, $1\le j \le N$, are independent, identically
distributed $\mathcal N_{\bC}(0, \sigma_{G_1}^2(r_0))$ with
\begin{equation}\label{eq18}
\sigma_{G_1}^2(r_0) = N\, \sum_{k\ge 1} a_{kN}^2 r_0^{2kN}\,.
\end{equation}
In addition, we have
\begin{equation}\label{eq19}
0 \le \sigma_F^2(r_0) - \sigma_{G_1}^2(r_0) \le C
\sigma_F^{2-c}(r_0)\,.
\end{equation}
\end{lemma}

\noindent{\em Proof of Lemma~\ref{lemma18a}}: Applying
Lemma~\ref{lemma14} to the function $G_1$ evaluated at $\tau z$ we see that the eigenvalues of the covariance
matrix of the random variables $(G_1(\omega^j\tau r_0))_{1\le j \le
N}$ are all equal to $N\, \sum_{k\ge 1} a_{kN}^2 r_0^{2kN}$. Hence, the covariance
matrix of these Gaussian random variables is diagonal, that is, they
are independent, and the relation \eqref{eq18} holds.

To prove estimate~\eqref{eq19}, observe that
\[
\sigma_{G_1}^2(r_0) \le \sum_{n\ge 0} a_n^2 r_0^{2n} = \sigma_{F}^2(r_0),
\]
and since the sequence $(a_n)$ does not increase,
we have
\[
\sigma_{G_1}^2(r_0) \ge \sum_{n\ge N} a_n^2 r_0^{2n} = \sigma_F^2(r_0) - \sum_{n=0}^{N-1} a_n^2 r_0^{2n}\,.
\]
Recalling (see \eqref{eq:a_n_asymptotics}) that
\[
a_n^2 = \frac{\Gamma(n+L)}{\Gamma(L)\Gamma(n+1)} \sim \frac{n^{L-1}}{\Gamma(L)},
\qquad n\to\infty\,,
\]
we see that
\[
\sum_{n=0}^{N-1} a_n^2 r_0^{2n} \le C\, \sum_{n=0}^{N-1} (1+n)^{L-1} \le C N^L \le C \sigma_F^{2-c}(r_0)\,,
\]
proving the lemma. \hfill $\Box$

\medskip
We henceforth condition on $F(0)$ and $G_2$ (that is, on $\zeta_0$ and $(\zeta_n'')_{1\le n \le N-1}$),
and write
\begin{align*}
\bE^{F(0), G_2} \bigl[\ . \ \bigr] &= \bE \bigl[\ . \ \big|\, F(0), G_2 \bigr], \\
\bP^{F(0), G_2} \bigl[\ . \ \bigr] &= \bP \bigl[\ . \ \big|\, F(0), G_2 \bigr].
\end{align*}

In the following section, we consider the case when $|F(0)|/\sigma_F$ is sufficiently small.

\subsection{$|F(0)|\le a \sigma_F$}\label{sec:F_0_small}
We show that the intersection of the hole event with the event $\{|F(0)|\le a \sigma_F\}$ is negligible for $a$ sufficiently small.

By Lemma~\ref{lemma9} (with $\kappa=2$ and $k=N$),
it suffices to estimate the probability
\[
\bP^{F(0), G_2} \Bigl[ \sum_{j=1}^N \log\Bigl| 1+\frac{G(\omega^j\tau r_0)}{F(0)} \Bigr| \le C \Bigr]
= \bP^{F(0), G_2} \Bigl[ \prod_{j=1}^N \Bigl| 1+\frac{G(\omega^j\tau r_0)}{F(0)} \Bigr|^{-1} \ge e^{-C} \Bigr]
\]
with some fixed $\tau\in \bT$. By Chebyshev's inequality, the right-hand side is bounded by
\[
C\,\bE^{F(0), G_2} \Bigl[ \prod_{j=1}^N \Bigl| 1+\frac{G(\omega^j\tau r_0)}{F(0)} \Bigr|^{-1} \Bigr]
= C\, \prod_{j=1}^N \bE^{F(0), G_2} \Bigl[ \Bigl| 1+\frac{G(\omega^j\tau r_0)}{F(0)} \Bigr|^{-1} \Bigr]
\]
where in the last equality we used the independence of $\bigl( G_1(\omega^j \tau r_0) \bigr)_{1\le j \le N}$ proven in
Lemma~\ref{lemma18a}, and the independence of $G_1$ and $F(0)$, $G_2$. Then, applying Lemma~\ref{lemma15} with
$t=|F(0)|/\sigma_{G_1}(r_0)$, we get using \eqref{eq19}, for each $1\le j \le N$,
\begin{multline*}
\bE^{F(0), G_2} \Bigl[ \Bigl| 1+\frac{G(\omega^j\tau r_0)}{F(0)} \Bigr|^{-1} \Bigr] \\
= \bE^{F(0), G_2} \Bigl[ \Bigl| 1+\frac{G_2(\omega^j\tau r_0)}{F(0)} + \frac{G_1(\omega^j\tau r_0)}{F(0)} \Bigr|^{-1} \Bigr]
\le C\, \frac{|F(0)|}{\sigma_{G_1}} \le C \cdot  a < \frac12\,,
\end{multline*}
provided that the constant $a$ is sufficiently small.

Thus,
\[
\bP^{F(0), G_2} \Bigl[ \prod_{j=1}^N \Bigl| 1+\frac{G(\omega^j\tau r_0)}{F(0)} \Bigr|^{-1} \ge e^{-C} \Bigr]
\le C 2^{-N} \le Ce^{-c\delta^{-\alpha}}\,.
\]
Since $\alpha>L$, this case gives a negligible contribution to the probability of the hole event.

\subsection{$a\sigma_F \le |F(0)| \le A\sigma_F\sqrt{\log\tfrac{1}{(1-r)\sigma_F^2}}$}
Our strategy is to make the constant $A$ as large as possible, while keeping the hole event negligible.
Then, up to negligible terms, we will bound $\bP[\text{Hole}(r)]$  by
\[
\bP\Big[ |F(0)|>A\sigma_F\sqrt{\log\tfrac{1}{(1-r)\sigma_F^2}} \Big] = \exp \Big[-A^2\sigma_F^2 \log\tfrac{1}{(1-r)\sigma_F^2} \Big]\,.
\]

We fix a small positive parameter $\e$, put
\[
J = \Bigl\{1\le j \le N\colon \Bigl| 1+\frac{G_2(\omega^j r_0)}{F(0)} \Bigr|\ge 1+2\e \Bigr\}\,,
\]
and introduce the event $\bigl\{ |J|\le (1-2\eta) N \bigr\}$,  where $\eta=\delta^{\al_0}$, $\al_0$ is a sufficiently small positive constant, and $|J|$ denotes the size of the set $J$.
We will estimate separately the probabilities of the events
\[
\text{Hole}(r)\bigcap \bigl\{ a\sigma_F \le |F(0)| \le A\sigma_F\sqrt{\log\tfrac{1}{(1-r)\sigma_F^2}} \bigr\}  \bigcap \bigl\{ |J|\le (1-2\eta) N \bigr\}
\]
and
\[
\text{Hole}(r)\bigcap \bigl\{ a\sigma_F \le |F(0)| \le A\sigma_F\sqrt{\log\tfrac{1}{(1-r)\sigma_F^2}} \bigr\} \bigcap \bigl\{ |J|> (1-2\eta) N \bigr\}\,.
\]

\medskip
Given $\tau\in\bT$, put
\begin{align*}
J_-(\tau) &= \Bigl\{1\le j \le N\colon \Bigl| 1+\frac{G_2(\omega^j \tau r_0)}{F(0)} \Bigr|\ge 1+\e \Bigr\}, \\
J_+(\tau) &= \Bigl\{1\le j \le N\colon \Bigl| 1+\frac{G_2(\omega^j \tau r_0)}{F(0)} \Bigr|\ge 1+3\e \Bigr\}.
\end{align*}
Our first goal is to show that, outside a negligible event, the random sets $J_\pm(\tau)$ are similar
in size to the set $J$.

\subsubsection{Controlling $G_2$ at the points $\bigl( \omega^j\tau r_0 \bigr)_{1\le j \le N}$}

\begin{lemma}\label{lemma3.1}
Given $\alpha_0>0$,
there is an event $\widetilde{\Omega}=\widetilde{\Omega}(r)$ with
$\bP^{F(0)} \bigl[ \widetilde{\Omega} \bigr] \le e^{-\delta^{-(L+c)}}$ on $\{|F(0)|\ge a\sigma_F\}$,
for $r$ sufficiently close to $1$, such that, for any
$\tau\in\bT$, on the complement $\widetilde{\Omega}^c$ we have
\begin{equation}\label{eq:J_-_+_events}
\begin{split}
  |J_-(\tau)|\ge |J| - \eta N,\\
  |J_+(\tau)|\le |J| + \eta N
\end{split}
\end{equation}
with $\eta=\delta^{\alpha_0}$.
\end{lemma}

\noindent{\em Proof of Lemma~\ref{lemma3.1}}: We take
small positive constants $\alpha_1$ and $\alpha_2$ satisfying
\begin{equation}\label{eq:alpha_assumptions}
  \begin{split}
    &\alpha_1<\max\left\{\frac{(1-\alpha)L}{2}, \frac{(1-L)\alpha_2}{2}\right\}\quad\text{and}\\
    &\alpha_1 + 2\alpha_2 < \alpha,
  \end{split}
\end{equation}
let $M=[\delta^{-\al_2}]$, write $G_2~=~G_3 + G_4$ as follows
\begin{align*}
G_3(z) &=\sum_{n=1}^M \zeta_n'' d_n z^n, \\
G_4(z) &=\sum_{n=M+1}^{N-1} \zeta_n'' d_n z^n,
\end{align*}
and define the events
\begin{align*}
\Omega_3 &= \bigcup_{n=1}^M
\bigl\{ |\zeta_n''| \ge \delta^{-\al_1} |F(0)| \bigr\}, \\
\Omega_4 &= \Bigl\{ \sum_{n=M+1}^{N-1} |\zeta_n''|^2 d_n^2 \ge \delta^{2 \al_1} |F(0)|^2\Bigr\}\,.
\end{align*}
First, we show that (under our running assumption $|F(0)|\ge a\sigma_F$)
the events $\Omega_3$ and $\Omega_4$ are negligible. Then, outside these events, we
estimate the functions $G_3$ and $G_4$ at the points $\bigl( \omega^j\tau r_0 \bigr)_{1\le j \le N}$.
We may assume without loss of generality that
$|\arg(\tau)|\le \pi/N$, as rotating $\tau$ by $2\pi k/N$ leaves $|J_\pm (\tau)|$ unchanged.

\begin{lemma}\label{lemma3.1a}
For $r$ sufficiently close to $1$, we have $\bP^{F(0)}\bigl[ \Omega_3 \bigr] \le e^{-\delta^{-(L+c)}}$  on $\{|F(0)|\ge a\sigma_F\}$.
\end{lemma}

\noindent{\em Proof of Lemma~\ref{lemma3.1a}}:
Using the union bound, we get
\[
\bP^{F(0)}\bigl[ \Omega_3 \bigr] \le
\sum_{n=1}^M \bP^{F(0)}\bigl[ |\zeta_n''|\ge \delta^{-\al_1} |F(0)| \bigr]
\le \sum_{n=1}^M \bP^{F(0)}\bigl[ |\zeta_n''|\ge \delta^{-\al_1} a \sigma_F \bigr]
\le Me^{-c a^2\delta^{-(L+2 \al_1)}}\,.\hfill \qed
\]

\begin{lemma}\label{lemma3.1b}
For $r$ sufficiently close to $1$, we have $\bP^{F(0)}\bigl[ \Omega_4 \bigr] \le e^{-\delta^{-(L+c)}}$ on $\{|F(0)|\ge a\sigma_F\}$.
\end{lemma}

\noindent{\em Proof of Lemma~\ref{lemma3.1b}}:
Put
\[
X=\sum_{n=M+1}^{N-1} |\zeta_n''|^2 d_n^2\,.
\]
Let $\la>0$ satisfy $\la d_n^2<1$ for $M+1\le n\le N-1$. Then,
\begin{multline*}
\bP \bigl[ X\ge t \bigr] = \bP \bigl[ e^{\la X} \ge e^{\la t} \bigr]
\le e^{-\la t} \bE \Bigl[ \prod_{n=M+1}^{N-1} e^{\la |\zeta_n''|^2 d_n^2} \Bigr] \\
= e^{-\la t} \prod_{n=M+1}^{N-1} \frac1{1-\la d_n^2}
= \exp \Bigl[ -\la t + \sum_{n=M+1}^{N-1} \log\frac1{1-\la d_n^2}\Bigr]\,.
\end{multline*}
We take $\la = 1/(2 a_M^2)$, recalling that $a_n^2 = b_n^2 + d_n^2$ with $(a_n)$ non-increasing. Then,
\[
\log\frac1{1-\la d_n^2} \le 2\la d_n^2\,.
\]
Thus,
\begin{equation}\label{eqA}
\bP \bigl[ X\ge t \bigr] \le
\exp \Bigl[ -\la \Bigl( t - 2\sum_{n=M+1}^{N-1} d_n^2 \Bigr) \Bigr]
\le \exp \bigl[ - \tfrac12 \la t\bigr]\,,
\end{equation}
provided that
\begin{equation}\label{eqB}
t\ge 4 \sum_{n=M+1}^{N-1} d_n^2\,.
\end{equation}
We use estimate~\eqref{eqA} with $t=\delta^{2 \al_1}|F(0)|^2$. Note that
$ \delta^{2 \al_1}|F(0)|^2 \ge c\delta^{-(L-2 \al_1)} $
and that
\[
4 \sum_{n=M+1}^{N-1} d_n^2 \le  4 \sum_{n=M+1}^{N-1} a_n^2 \le C N^L < C \delta^{-\alpha L}\,,
\]
which satisfies \eqref{eqB} by \eqref{eq:alpha_assumptions}.
Finally, using again \eqref{eq:alpha_assumptions}, we get
\begin{multline*}
\bP \bigl[ X \ge \delta^{2 \al_1} |F(0)|^2 \bigr]
\le \exp \bigl[ -\tfrac1{(2a_M)^2} \cdot \delta^{2 \al_1} |F(0)|^2 \bigr] \\ \\
\le \exp \bigl[ -c  M^{1-L} \delta^{2 \al_1} |F(0)|^2 \bigr]
\le \exp \bigl[ -ca^2 \delta^{-L-(1-L) \al_2 +2 \al_1} \bigr] \le
\exp\bigl[-\delta^{-(L+c)}\bigr].\qed
\end{multline*}

\medskip

\begin{lemma}\label{lemma3.1c} Suppose that $r$ is sufficiently close to $1$. Then on the event
$\Omega_3^c$ we have
\[
\sup_{|\arg \tau|<\pi N^{-1}}\, \sup_{1\le j \le N} \Bigl| G_3(\omega^j\tau r_0) - G_3(\omega^jr_0) \Bigr| < \frac12\, \e\,|F(0)|,
\qquad 1 \le j \le N.
\]
\end{lemma}

\noindent{\em Proof of Lemma~\ref{lemma3.1c}}: For each $1\le j \le N$, since $0 \le d_n \le 1$, we have
\begin{multline*}
\bigl|  G_3(\omega^j\tau r_0)- G_3(\omega^jr_0) \bigr| \le
\max_{\bar\bD} |G_3'| \cdot \frac{\pi}N \le
M\, \sum_{n=1}^M |\zeta_n''|\cdot |d_n| \cdot \frac{\pi}N \\
\le \frac{CM^2}{N}\cdot \delta^{- \al_1} |F(0)| \le C\delta^{-2 \al_2+\alpha- \al_1} |F(0)|
<\tfrac12 \e |F(0)|\,,
\end{multline*}
provided that $2 \al_2 + \al_1 < \alpha$, and that $r$ is sufficiently close to $1$.
\hfill $\Box$

\begin{lemma}\label{lemma3.1d}
Suppose that $r$ is sufficiently close to $1$.
Then, on the event $\Omega_4^c$, for any $\tau\in\bT$,  the cardinality of the set
\[
\Bigl\{ 1\le j \le N\colon \max\bigl( |G_4(\omega^jr_0)|, |G_4(\omega^j\tau r_0)| \bigr) \ge
\tfrac14\, \e |F(0)| \Bigr\}
\]
does not exceed $\eta N$.
\end{lemma}
\noindent{\em Proof of Lemma~\ref{lemma3.1d}}:
We have
\begin{align*}
\frac1{N}\, \sum_{j=1}^{N} &|G_4(\omega^j r_0 )|^2
= \frac1{N}\, \sum_{j=1}^{N} \biggl| \sum_{n=M+1}^{N-1} \zeta_n'' d_n (\omega^j r_0)^n  \biggr|^2 \\
&=\sum_{n_1, n_2 = M+1}^{N-1} \zeta_{n_1}'' \overline{\zeta_{n_2}''} d_{n_1} d_{n_2} r_0^{n_1+n_2}\,
\frac1{N}\,\sum_{j=1}^N \omega^{j(n_1 - n_2)} \\
&= \sum_{n=M+1}^{N-1} |\zeta_n''|^2 d_n^2 r_0^{2n} \\
&\le \sum_{n=M+1}^{N-1} |\zeta_n''|^2 d_n^2 \le \delta^{2 \al_1} |F(0)|^2\qquad\text{on }\Omega_4^c\,,
\end{align*}
and similarly,
\[
\frac1{N}\, \sum_{j=1}^{N} |G_4(\omega^j \tau r_0 )|^2 \le \delta^{2 \al_1} |F(0)|^2\qquad\text{on }\Omega_4^c\,.
\]
Hence, on $\Omega_4^c$, the cardinality of the set we are interested in does not exceed
\[
\frac{32\,\delta^{2 \al_1}N}{\e^2} < \delta^{\al_0}N = \eta N\,,
\]
provided that $2 \al_1<\al_0$ and that $r$ is sufficiently close to $1$.
\hfill $\Box$

\medskip Now, Lemma~\ref{lemma3.1} is a straightforward consequence of Lemmas~\ref{lemma3.1a},
~\ref{lemma3.1b}, ~\ref{lemma3.1c}, and~\ref{lemma3.1d}. \hfill $\Box$

\subsubsection{$|J|> (1-2\eta)N$, $\eta=\delta^{\alpha_0}$}
In this section we show that the intersection of the hole event with the event
\begin{equation}\label{eq:F_J_event}
  \Big\{a\sigma_F \le |F(0)| \le A\sigma_F\sqrt{\log\tfrac{1}{(1-r)\sigma_F^2}}\Big\}\cap\Big\{|J|> (1-2\eta)N\Big\}
\end{equation}
is negligible. Taking into account the fact that $J, J_-(\tau)$ and $J_+(\tau)$ are measurable with respect to $F(0)$ and $G_2$, Lemma~\ref{lemma9} (with $\kappa=2$ and $k=N$) and
Lemma~\ref{lemma3.1} show that it suffices to estimate uniformly in $\tau\in\bT$, on the intersection of the events \eqref{eq:F_J_event} and \eqref{eq:J_-_+_events}, the
probability
\begin{equation*}
\bP^{F(0), G_2} \Bigl[ \sum_{j=1}^N \log\Bigl| 1+\frac{G(\omega^j\tau r_0)}{F(0)} \Bigr| \le C \Bigr].
\end{equation*}
Taking some positive $\theta_1=\theta_1(r)$ tending to $0$ as $r\to 1$ (the function $\theta_1(r)$ will be chosen later) and applying Chebyshev's inequality, the last probability does
not exceed
\begin{multline}\label{eqC}
C\,\bE^{F(0), G_2} \Bigl[ \prod_{j=1}^N \Bigl| 1+\frac{G(\omega^j\tau r_0)}{F(0)} \Bigr|^{-\theta_1} \Bigr]
= C\,\bE^{F(0), G_2} \Bigl[ \prod_{j=1}^N \Bigl| 1+\frac{G_2(\omega^j\tau r_0)}{F(0)}
+\frac{G_1(\omega^j\tau r_0)}{F(0)}\Bigr|^{-\theta_1} \Bigr] \\
= C\, \prod_{j=1}^N \bE^{F(0), G_2} \Bigl[ \Bigl| 1+\frac{G_2(\omega^j\tau r_0)}{F(0)}
+\frac{G_1(\omega^j\tau r_0)}{F(0)} \Bigr|^{-\theta_1} \Bigr]\,.
\end{multline}
Once again, we used the independence of $\bigl( G_1(\omega^j \tau r_0) \bigr)_{1\le j \le N}$
and the independence of $G_1$ and $F(0)$, $G_2$.

\medskip
For $j\in J_-(\tau)$, we have
\begin{multline*}
 \Bigl| 1+\frac{G_2(\omega^j\tau r_0)}{F(0)}
+\frac{G_1(\omega^j\tau r_0)}{F(0)}  \Bigr| =
\Bigl| 1+\frac{G_2(\omega^j\tau r_0)}{F(0)} \Bigr| \cdot
\Bigl| 1+\frac{G_1(\omega^j\tau r_0)}{F(0)+G_2(\omega^j \tau r_0)} \Bigr| \\
\ge (1+\e) \Bigl| 1+\frac{G_1(\omega^j\tau r_0)}{F(0)+G_2(\omega^j \tau r_0)} \Bigr|\,.
\end{multline*}
Hence, by Lemma~\ref{lemma17}, for such $j$ we obtain
\begin{equation*}
\bE^{F(0), G_2} \Bigl[ \Bigl| 1+\frac{G_2(\omega^j\tau r_0)}{F(0)}
+\frac{G_1(\omega^j\tau r_0)}{F(0)} \Bigr|^{-\theta_1} \Bigr]
\le (1+\e)^{-\theta_1} (1+C\theta_1^2) \le e^{-c\e\theta_1 +C\theta_1^2}
< e^{-c\e\theta_1},
\end{equation*}
provided that $r$ is so close to $1$ that
$ \theta_1 (r)$ is much smaller than $\e$ (recall that $\e$ is small but fixed).

\medskip
For $j\notin J_-(\tau)$, using Lemma~\ref{lemma15} (with $t = \frac{|F(0)|}{\sigma_G}
\le CA\,\sqrt{\log\tfrac{1}{(1-r)\sigma_F^2}}$), we get
\[
\bE^{F(0), G_2} \Bigl[ \Bigl| 1+\frac{G(\omega^j\tau r_0)}{F(0)} \Bigr|^{-\theta_1} \Bigr]
\le (CA\,\sqrt{\log\tfrac{1}{(1-r)\sigma_F^2}})^{\theta_1} (1+C\theta_1) < C( \log\tfrac{1}{(1-r)\sigma_F^2} )^{\frac12 \theta_1}\,.
\]
Thus, \eqref{eqC} does not exceed
\[
\exp\bigl[ -c\e\theta_1 |J_-(\tau)| + \Big(C + \frac12 \theta_1 \log\log\tfrac{1}{(1-r)\sigma_F^2} \Big) (N-|J_-(\tau)|) \bigr]\,.
\]
As we are on the intersection of the events \eqref{eq:F_J_event} and \eqref{eq:J_-_+_events} we have
$|J_-(\tau)| \ge (1-3\eta)N$. Therefore, the expression in the last displayed formula
does not exceed
\[
\exp\bigl[ -c\e\theta_1 N +C \eta N (1 + \theta_1 \log\log\tfrac{1}{(1-r)\sigma_F^2}) \bigr] \stackrel{\rm def}= E.
\]
Then, letting $\theta_1 = \delta^c$ with $c<\min\{\al_0, \alpha - L\}$, and using that
\[
\theta_1 \log\log\tfrac{1}{(1-r)\sigma_F^2}  \le \delta^c \bigl( \log\log \frac1{\delta} + C \bigr)\to 0 \quad \text{as}\
\delta\to 0\,,
\]
we see that $E \le \exp\bigl[ -c\e\theta_1 N \bigr]$ (recall that $\eta = \delta^{\al_0}$) and conclude that the event
\[
\text{Hole}(r) \bigcap \bigl\{ a\sigma_F \le |F(0)| \le A\sigma_F\sqrt{\log\tfrac{1}{(1-r)\sigma_F^2}}\bigr\} \bigcap\bigl\{ |J|> (1-2\eta) N \bigr\}
\]
is negligible.

\subsubsection{$|J| \le (1-2\eta)N$, $\eta=\delta^{\alpha_0}$}\label{subsubsect:smallJ}

Here we show that the intersection of the hole event with the event
\begin{equation}\label{eq:F_J_event2}
  \Big\{a\sigma_F \le |F(0)| \le A\sigma_F\sqrt{\log\tfrac{1}{(1-r)\sigma_F^2}} \Big\}\cap\Big\{|J|\le (1-2\eta)N\Big\}
\end{equation}
is negligible, provided that $ A^2 < \frac{1}{2}$.

In this case, our starting point is Lemma~\ref{lemma8} which we apply with
the polynomial $P=F(0)+G_2$. Combined with Lemma~\ref{lemma5} and Lemma~\ref{lemma3.1}, it tells us
that it suffices to estimate uniformly in $\tau\in\bT$, on the intersection of the events \eqref{eq:F_J_event2} and \eqref{eq:J_-_+_events}, the probability
\[
\bP^{F(0), G_2} \Bigl[ \sum_{j=1}^N \log \Bigl| \frac{F}{P} (\omega^j\tau r_0)\Bigr| \le C \Bigr]\,.
\]
Noting that
\[
\frac{F}P = 1+ \frac{G_1}{F(0)+G_2},
\]
we rewrite this expression as
\begin{align*}
\bP^{F(0), G_2} \Bigl[ \prod_{j=1}^N \Bigl|  1 &+ \frac{G_1(\omega^j\tau r_0)}{F(0)+G_2(\omega^j\tau r_0)} \Bigr|
\le e^{C} \Bigr] \\
&= \bP^{F(0), G_2} \Bigl[ \prod_{j=1}^N \Bigl|  1+ \frac{G_1(\omega^j\tau r_0)}{F(0)+G_2(\omega^j\tau r_0)} \Bigr|^{-\theta_2}
\ge e^{-C\theta_2} \Bigr] \\
&\le e^{C\theta_2}\, \bE^{F(0), G_2} \Bigl[ \prod_{j=1}^N \Bigl|  1+ \frac{G_1(\omega^j\tau r_0)}{F(0)+G_2(\omega^j\tau r_0)} \Bigr|^{-\theta_2} \Bigr]\,.
\end{align*}
The positive $\theta_2=\theta_2(r)$ (again tending to $0$ as $r\to 1$) will be chosen later.
By the independence of $(G_1(\omega^j \tau r_0)_{1\le j \le N}$ proven in
Lemma~\ref{lemma18a}, and the independence of $G_1$ and $F(0)$, $G_2$, it suffices to
estimate the product
\[
\prod_{j=1}^N  \bE^{F(0), G_2} \Bigl[ \Bigl|  1+ \frac{G_1(\omega^j\tau r_0)}{F(0)+G_2(\omega^j\tau r_0)} \Bigr|^{-\theta_2} \Bigr]\,.
\]

\medskip
First we consider the terms with $j\in J_+(\tau)$. In this case, by Lemma~\ref{lemma17}, we have
\begin{equation}\label{eqG1}
\bE^{F(0), G_2} \Bigl[ \Bigl|  1+ \frac{G_1(\omega^j\tau r_0)}{F(0)+G_2(\omega^j\tau r_0)} \Bigr|^{-\theta_2} \Bigr]
\le 1+ C\theta_2^2 < e^{C\theta_2^2}\,.
\end{equation}
For the terms with $j\notin J_+(\tau)$, we apply Lemma~\ref{lemma17} with
\begin{multline*}
\qquad t = \frac{|F(0)+G_2(\omega^j\tau r_0)|}{\sigma_{G_1}} \le (1+3\e)  \frac{|F(0)|}{\sigma_{G_1}} \\
= (1+3\e) \frac{\sigma_F}{\sigma_{G_1}}\, \frac{|F(0)|}{\sigma_F}
\stackrel{\eqref{eq19}}\le  (1+4\e) A \sqrt{\log\tfrac{1}{(1-r)\sigma_F^2}}\,, \qquad
\end{multline*}
provided that $r$ is sufficiently close to $1$. Then, by Lemma~\ref{lemma17} (using that $(1-r)\sigma_F^2$ is sufficiently small)
\[
\bE^{F(0), G_2} \Bigl[ \Bigl|  1+ \frac{G_1(\omega^j\tau r_0)}{F(0)+G_2(\omega^j\tau r_0)} \Bigr|^{-\theta_2} \Bigr]
\le 1 - c\theta_2 \big((1-r)\sigma_F^2\big)^{(1+10\e) A^2} + C\theta_2^2\,.
\]
Assuming that
\begin{equation}\label{eqH}
\theta_2 =o(1) \big((1-r)\sigma_F^2\big)^{(1+10\e) A^2}\,,
\end{equation}
we continue our estimate as follows
\begin{equation}\label{eqJ}
\le 1 - c\theta_2 \big((1-r)\sigma_F^2\big)^{(1+10\e) A^2} \le \exp\bigl[ -c\theta_2 \big((1-r)\sigma_F^2\big)^{(1+10\e) A^2} \bigr]\,.
\end{equation}

Multiplying the bounds~\eqref{eqG1} and~\eqref{eqJ}, we get
\begin{multline*}
\prod_{j=1}^N  \bE^{F(0), G_2} \Bigl[ \Bigl|  1+ \frac{G_1(\omega^j\tau r_0)}{F(0)+G_2(\omega^j\tau r_0)} \Bigr|^{-\theta_2} \Bigr] \\
\le \exp \bigl[ C\theta_2^2 |J_+(\tau)| - c\theta_2 \big((1-r)\sigma_F^2\big)^{(1+10\e) A^2} (N-|J_+(\tau)|) \bigr]\,.
\end{multline*}
As we are on the intersection of the events \eqref{eq:F_J_event2} and \eqref{eq:J_-_+_events}, we have $|J_+(\tau)| \le |J| + \eta N \le (1-\eta)N$. Then, the expression in the exponent
on the right-hand side of the previous displayed equation does not exceed
\[
\bigl( C\theta_2^2  - c\eta \theta_2 \big((1-r)\sigma_F^2\big)^{(1+10\e) A^2} \bigr) N\,.
\]
Letting
$ \theta_2 = c \eta \big((1-r)\sigma_F^2\big)^{(1+10\e) A^2} $
(which satisfies our previous requirement~\eqref{eqH} since $\eta\to 0$ as $r\to 1$), we estimate the previous expression
by
\[
-c \eta^2 \big((1-r)\sigma_F^2\big)^{2(1+10\e) A^2} N = -c \delta^{2 \al_0+2(1+10\e) A^2 (1-L) -\alpha}\,.
\]
To make the event with probability bounded by $\exp\bigl[ -c \delta^{2 \al_0+2(1+10\e) A^2 (1-L) -\alpha} \bigr]$ negligible,
we need to be sure that $\alpha - 2(1+10\e) A^2(1-L) - 2 \al_0 > L$. Since the constants $\e$ and $\al_0$ can be made arbitrarily small,
while $\alpha$ can be made arbitrarily close to $1$, we conclude that the event
\[
\text{Hole}(r)\bigcap \bigl\{ a\sigma_F \le |F(0)| \le A\sigma_F\sqrt{\log\tfrac{1}{(1-r)\sigma_F^2}}\bigr\} \bigcap \{|J| \le (1-2\eta)N\}
\]
is negligible, provided  that $ A^2 < \frac{1}{2}$.

\medskip
We conclude that the event
\[
\text{Hole}(r)\bigcap \bigl\{ a \sigma_F <|F(0)|<A\sigma_F\sqrt{\log\tfrac{1}{(1-r)\sigma_F^2}}\bigr\}
\]
is negligible whenever $ A^2 < \frac{1}{2}$, and therefore, combined with the bound of Section~\ref{sec:F_0_small},
for any such $A$,
\begin{align*}
\bP \bigl[ \text{Hole}(r) \bigr] & \le \bP\Bigl[ |F(0)|\ge A\sigma_F\sqrt{\log\tfrac{1}{(1-r)\sigma_F^2}}\, \Bigr] + \langle \text{negligible\ terms} \rangle \\
& = \exp\big(-A^2\sigma_F^2 \log \tfrac{1}{(1-r)\sigma_F^2}\big) + \langle \text{negligible\ terms} \rangle \,,
\end{align*}
whence, letting $r\uparrow 1$,
\begin{align*}
\bP \bigl[ \text{Hole}(r) \bigr]
&\le \exp\Bigl[ - \frac{1-o(1)}{2}\, \sigma_F^2 \log\tfrac{1}{(1-r)\sigma_F^2}  \Bigr]  \\
&= \exp\Bigl[ - \frac{1-o(1)}{2} \cdot \frac{1}{(2\delta)^{L}} \cdot  (1-L) \log\frac1{\delta}  \Bigr] \\
&= \exp\Bigl[ - \frac{1-L-o(1)}{2^{L+1}} \cdot \frac{1}{\delta^{L}} \, \log\frac1{\delta}  \Bigr]
\end{align*}
completing the proof of the upper bound in the case $0<L<1$.  \hfill $\Box$

\section{Lower bound on the hole probability for $0<L<1$}\label{section:LB,L<1}
As before, let $F=F(0)+G$. Fix $\e>0$, set
\[
M = \frac{\sqrt{1-L+2\e}}{(1-r^2)^{L/2}} \cdot \sqrt{\log\frac1{1-r}}\,,
\]
define the events
\[
\Omega_1 = \bigl\{ |F(0)| > M \bigr\}, \qquad
\Omega_2 = \bigl\{ \max_{r\bar\bD} |G| \le M  \bigr\} = \bigl\{ \max_{r \bT} |G| \le M  \bigr\}\,,
\]
and observe that $\text{Hole}(r) \supset \Omega_1 \bigcap \Omega_2$ and that $\Omega_1$ and $\Omega_2$ are independent.

Put
\[
N = \bigl[ (1-r)^{-1-\e} \bigr]\,, \quad \omega = e(1/N)\,, \quad z_j = r\omega^j, \ 1\le j \le N\,,
\]
\[
M' = \frac{\sqrt{1-L+\e}}{(1-r^2)^{L/2}}\, \sqrt{\log\frac1{1-r}}\,, \quad
M'' = \frac{A}{(1-r^2)^{(L+2)/2}}\, \sqrt{\log\frac1{1-r}}\,,
\]
with a sufficiently large positive constant $A$ that will be chosen later,
and define the events
\[
\Omega_3 = \Bigl\{ \max_{1\le j \le N } |G(z_j)| \le M'  \Bigr\}\,, \quad
\Omega_4 = \Bigl\{ \max_{r\bT} |G'| \le M''  \Bigr\}\,.
\]
Then, $ \Omega_2 \supset \Omega_3 \bigcap \Omega_4$ when $r$ is sufficiently close to $1$ as a function of $\e$.
Hence,
\[
\bP \bigl[ \text{Hole}(r) \bigr] \ge \bP \bigl[ \Omega_1 \bigr] \cdot \bP\bigl[ \Omega_3 \bigcap \Omega_4\bigr]\,.
\]
For $K \in \bN$, we put $\la = e\big(1/2^K\big)$, $w_k = \la^k r$, $1 \le k \le 2^K$ and define
\[
\Omega_4^K = \Bigl\{ \max_{1\le k \le 2^K } |G'(w_k)| \le M''  \Bigr\}.
\]
Notice that $\Omega_4^K \downarrow \Omega_4$ as $K\to\infty$.
In order to bound $\bP\bigl[\, \Omega_3 \bigcap \Omega_4 \, \bigr]$ from below we use Harg\'e's version~\cite{Harge} of the
Gaussian correlation inequality\footnote{
Recently, Royen~\cite{Royen} proved the full version of the Gaussian correlation inequality, see also
the paper by Lata{\l}a and Matlak~\cite{LM}.
}:

\bigskip\noindent{\bf Theorem} (G. Harg\'e) {\em Let $\gamma$ be a Gaussian measure on $\bR^n$, let
$A\subset \bR^n$ be a convex symmetric set, and let $B\subset\bR^n$ be an ellipsoid (that is, a set
of the form $\{X\in\bR^n\colon \langle CX, X \rangle \le 1\}$, where $C$ is a non-negative symmetric matrix).
Then $ \gamma (A\bigcap B) \ge \gamma (A) \gamma (B)$. }

\bigskip
We apply this inequality $N$ times to the Gaussian measure on $\bR^{2\big(N - j + 1 + 2^K\big)}$, $1 \le j \le N$, generated by
\begin{align*}
    X_j^r=\text{Re}\,G\left(z_j\right), & \dots, X_N^r=\text{Re}\,G\left(z_N\right),\\
    X_j^i=\text{Im}\,G\left(z_j\right), & \dots, X_N^i=\text{Im}\,G\left(z_N\right),
\end{align*}
and
\begin{align*}
    X_{N+1}^r=\text{Re}\,G'(w_1), & \dots, X_{N+2^{K}}^r=\text{Re}\,G'\left(w_{2^{K}}\right),\\
    X_{N+1}^i=\text{Im}\,G'(w_1), & \dots, X_{N+2^{K}}^i=\text{Im}\,G'\left(w_{2^{K}}\right),
\end{align*}
and the sets
\begin{align*}
A_j & = \left\{ \max_{j + 1\le k \le N} |G(z_j)| \le M^\prime \right\} \cap \left\{ \max_{1\le k \le 2^K} |G^\prime(w_k)| \le M^{\prime\prime} \right\},\\ \\
B_j & = \left\{ |G(z_j)|^2 = (X_j^r)^2 + (X_j^i)^2 \le (M^\prime)^2 \right\}.
\end{align*}
\noindent Thus, we get
\[
\bP [\Omega_3 \cap \Omega_4^K] \ge \prod_{j=1}^N \bP\bigl[ |G(z_j)| \le M'\, \bigr] \cdot \bP [\Omega_4^K]\,.
\]
Thus, by the monotone convergence of $\Omega_4^K$ to $\Omega_4$,
\[
\bP\bigl[ \Omega_3 \bigcap \Omega_4\bigr] \ge \bP \bigl[  \max_{r\bar\bD} |G'| \le M'' \bigr] \cdot
\prod_{j=1}^N \bP\bigl[ |G(z_j)| \le M' \bigr]\,.
\]

For each $1\le j \le N$, $G(z_j)$ is a complex Gaussian random variable with variance at most
$(1-r^2)^{-L}$, so that
\begin{multline*}
\prod_{j=1}^N \bP\bigl[ |G(z_j)| \le M' \bigr]
\ge \Bigl( 1-e^{-M'^2(1-r^2)^L} \Bigr)^N
= \Bigl( 1-e^{(1-L+\e)\,\log(1-r)} \Bigr)^N
= \Bigl( 1-(1-r)^{1-L+\e} \Bigr)^N \\
\ge  \exp\bigl[ -c N (1-r)^{1-L+\e} \bigr]
= \exp\bigl[ -c (1-r)^{-1-\e+1-L+\e}  \bigr]
= \exp\bigl[ -c (1-r)^{-L}  \bigr]
\end{multline*}
with $r$ sufficiently close to $1$.

Next note that
\[
G'(z) = \sum_{n\ge 0} \zeta_{n+1} (n+1)a_{n+1} z^n\,,
\]
and therefore we have
\[
\max_{s\bar\bD} \bE \bigl[ |G'|^2 \bigr] \le C (1-s^2)^{-(L+2)}\,, \qquad 0<s<1\,.
\]
Thus, noting that
\[
M'' \ge  cA \sqrt{\log (1-r)^{-1}} \cdot \sigma_{G'}(\tfrac12 (1+r))
\]
and applying Lemma~\ref{lemma2}, we see that
\[
\bP \bigl[ \max_{r\bar\bD} |G'| > M'' \bigr]
\le \frac{C}{1-r} \exp \bigl[ c A^2 \log (1-r) \bigr] < \tfrac12\,,
\]
provided that the constant $A$ in the definition of $M''$ is chosen sufficiently large.
Thus, $ \bP \bigl[ \max_{r\bar\bD} |G'| \le M'' \bigr] \ge \tfrac12$. Then, piecing everything together,
we get
\begin{align*}
\bP \bigl[ \text{Hole}(r) \bigr] &\ge \bP \bigl[ \Omega_1 \bigr] \cdot \bP\bigl[ \Omega_3 \bigcap \Omega_4\bigr] \\
&\ge  \bP \bigl[ \Omega_1 \bigr] \cdot
\bP \bigl[  \max_{r\bar\bD} |G'| \le M'' \bigr] \cdot
\prod_{j=1}^N \bP\bigl[ |G(z_j)| \le M' \bigr] \\
&\ge \bP \bigl[ \Omega_1 \bigr] \cdot \tfrac12\, \exp\bigl[ -c (1-r)^{-L}\bigr]\,.
\end{align*}
Finally, the theorem follows from the fact that
$  \bP \bigl[ \Omega_1 \bigr] = e^{-M^2} $. \hfill $\Box$

\subsection{Remark}\label{subsect:remark}
Having in mind the gap between the upper and lower bounds on the hole probability for $0<L<1$, as given in Theorem~\ref{thm:main}, we note here that a different method would be required
in order to improve the lower bound. Precisely, setting $F = F(0) + G$ as before, we will show that
\begin{multline}\label{eq:lower_bound_obstruction}
  \sup_M \bP\Bigl[|F(0)|>M,\;\; \max_{r\bar \bD} |G|\le M\Bigr]\\
   = \exp\left(-\frac{1-L+o(1)}{2^{L}}\, \frac1{(1-r)^L}\, \log\frac1{1-r}\right), \quad\text{as\ } r\uparrow 1.
\end{multline}
Our proof above shows that this holds with a greater or equal sign instead of the equality sign and it remains to establish the opposite inequality. It is clear that
\begin{equation*}
  \bP\Bigl[|F(0)|>M,\;\; \max_{r\bar \bD} |G|\le M\Bigr]\le \bP\Bigl[|F(0)|>M\Bigr] = e^{-M^2}
\end{equation*}
and hence the opposite inequality need only be verified in the regime where
\begin{equation*}
  M\le \frac{\sqrt{1-L-o(1)}}{(1-r^2)^{L/2}} \cdot \sqrt{\log\frac1{1-r}}, \quad\text{as\ } r\uparrow 1.
\end{equation*}
Now let $0<\e<\tfrac12 (1-L)$ and suppose that
\begin{equation}\label{eq:M_upper_bound}
  M\le \frac{\sqrt{1-L-2\e}}{(1-r^2)^{L/2}} \cdot \sqrt{\log\frac1{1-r}}.
\end{equation}
Set also $N=\bigl[ (1-r)^{-1+\epsilon}\bigr]$ and $\omega = e(1/N)$. Let $G_1$ and $G_2$ be as in Section~\ref{sec:splitting_G}, so that $G = G_1 + G_2$ and the random variables
$(G_1(\omega^j r))$, $1\le j\le N$, are independent and identically distributed by Lemma~\ref{lemma18a}. (The decomposition in Section~\ref{sec:splitting_G} yields independent random
variables at radius $r_0 = 1 - 2(1-r)$ but may easily be modified to radius $r$ (or indeed any radius).)
\[
  \bP\Bigl[|F(0)|>M,\;\; \max_{r\bar \bD} |G|\le M\Bigr]\le \bP\Bigl[\max_{1\le j\le N} |G_1(\omega^j r) +G_2(\omega^j r)|\le M\Bigr]\,.
\]
We condition on $G_2$ and write $\bP^{G_2}\bigl[\ . \ \bigr]=\bP \bigl[\ . \ \bigl|\, G_2\,\bigr]$.
Recalling that $(G_1(\omega^j r))$ are independent and applying the Hardy-Littlewood rearrangement inequality, we get
\begin{multline*}
\bP^{G_2}\Bigl[\max_{1\le j\le N} |G_1(\omega^j r) +G_2(\omega^j r)|\le M\Bigr]
= \prod_{j=1}^N \bP^{G_2}\Bigl[|G_1(\omega^j r) + G_2(\omega^j r)|\le M\, \Bigr] \\
\le \prod_{j=1}^N \bP^{G_2}\Bigl[|G_1(\omega^j r)|\le M\Bigr] = \left(1 - e^{-M^2 / \sigma_{G_1}^2(r)}\right)^N\le e^{-N\exp(-M^2 / \sigma_{G_1}^2(r))}\,.
\end{multline*}
Since the right-hand side does not depend on $G_2$, we can drop the conditioning on the left-hand side and finally get
\[
 \bP\Bigl[|F(0)|>M,\;\; \max_{r\bar \bD} |G|\le M\Bigr]\le e^{-N\exp(-M^2 / \sigma_{G_1}^2(r))}\,.
\]
It remains to note that with our choice of $N$ and using the upper bound \eqref{eq:M_upper_bound} on $M$ and the relation of $\sigma_F^2(r)$ and $\sigma_{G_1}^2(r)$ given in
Lemma~\ref{lemma18a} we have
\begin{equation*}
  N\exp(-M^2 / \sigma_{G_1}^2(r))\ge \Bigl(\frac{1}{1-r}\Bigr)^{L+\e-o(1)}, \quad\text{as\ } r\uparrow 1.
\end{equation*}
We conclude that
\begin{equation*}
  \bP\Bigl[|F(0)|>M,\;\; \max_{r\bar \bD} |G|\le M\Bigr] \le \exp\Bigl[ -\Bigl(\frac{1}{1-r}\Bigr)^{L+\e-o(1)} \Bigr], \quad\text{as\ } r\uparrow 1.
\end{equation*}
for all $M$ satisfying the upper bound \eqref{eq:M_upper_bound}. As $\e$ can be taken arbitrarily small, this completes the proof of \eqref{eq:lower_bound_obstruction}.

\section{Upper bound on the hole probability for $L>1$}\label{Sect:UB_L>1}

\subsection{Beginning the proof}
Put $\delta=1-r$, $1 < \kappa \le 2$ and $r_0=1-\kappa \delta$. We take
\begin{equation}\label{eq:choice_of_N}
N = \Bigl[\,\frac{L-1}{2\delta}\, \log\frac1{\delta}\,\Bigr]
\end{equation}
and put $\omega=e(1/N)$.
By Lemma~\ref{lemma9}, it suffices to estimate the probability
\[
\sup_{\tau\in\bT}\, \bP \Bigl[ \sum_{j=1}^N \log|F(\omega^j \tau r_0 )| \le N \log |F(0)| + C \Bigr]\,.
\]
Set
\begin{equation}\label{eq:choice_of_a}
a = \left(\log \frac1\delta \right)^{-\frac12}.
\end{equation}
Discarding a negligible event, we assume that
$ |F(0)| \le \delta^{-(\frac12 + a)} $. Let $0<\theta<2$ be a parameter depending on $\delta$ which we will
choose later. Then,
\begin{align*}
\bP \Bigl[ \sum_{j=1}^N \log|F(\omega^j \tau r_0 )| \le N \log |F(0)| + C \Bigr]
& \le \bP \Bigl[ \prod_{j=1}^N \bigl| F(\omega^j \tau r_0 ) \bigr|^{-\theta} \ge C \delta^{N\theta (\frac12+a)} \Bigr] \\
& \le C \delta^{-N\theta (\frac12+a)}\, \bE \Bigl[ \prod_{j=1}^N \bigl| F(\omega^j \tau r_0 ) \bigr|^{-\theta}\Bigr]\,.
\end{align*}
Using Lemma~\ref{lemma18}, we can bound the expectation on the right by
\[
\frac1{\det\Sigma}\, \Big( \Lambda^{\big(1-\frac12\theta\big)} \Gamma \bigl(1-\tfrac12\theta \bigr) \Big)^N \,,
\]
where $\Sigma$ is the covariance matrix of $\bigl( F(\omega^j \tau r_0) \bigr)_{1\le j \le N} $
and $\Lambda$ is the maximal eigenvalue of $\Sigma$. Note that, since the distribution of $F(z)$ is rotation
invariant, the covariance matrix $\Sigma$ does not depend on $\tau$.

\subsection{Estimating the eigenvalues of $\Sigma$}

By Lemma~\ref{lemma14}, the eigenvalues of $\Sigma$ are
\[
\lambda_m(\Sigma) = N\, \sum_{n\equiv m\, (N)} a_n^2 r_0^{2n}
= N \Bigl( a_m^2 r_0^{2m} + \sum_{j\ge 1} a_{m+jN}^2 r_0^{2(m+jN)} \Bigr)\,,
\qquad m=0, 1, \ldots , N-1\,.
\]
Now, for small $\delta$ and $j\ge 1$, we have
\begin{multline*}
\frac{a_{m+(j+1)N}^2 r_0^{2(m+(j+1)N)}}{a_{m+jN}^2 r_0^{2(m+jN)}} \le C
\Bigl( \frac{m+(j+1)N}{m+jN} \Bigr)^{L-1}\, r_0^{2N} \\
= C \Bigl( 1 + \frac{N}{m+jN} \Bigr)^{L-1}\, \bigl( 1-\kappa\delta\bigr)^{2N}
\le C e^{-2\kappa\delta N}\,.
\end{multline*}
Take $\delta$ sufficiently small so that
\[
\frac{a_{m+(j+1)N}^2 r_0^{2(m+(j+1)N)}}{a_{m+jN}^2 r_0^{2(m+jN)}} \le \frac12\,,
\]
which yields
\begin{align*}
\lambda_m(\Sigma) &\le N \Bigl( a_m^2 (1-\kappa\delta)^{2m}
+ a_{m+N}^2(1-\kappa\delta)^{2(m+N)} \sum_{j\ge 1} 2^{-j}\Bigr) \\
& = N \bigl( a_m^2 (1-\kappa\delta)^{2m}
+ a_{m+N}^2(1-\kappa\delta)^{2(m+N)} \bigr)\,.
\end{align*}

Put
\[
M \stackrel{\rm def}= \frac{Lr_0^2 -1}{1-r_0^2} = (1+o(1))\, \frac{L-1}{2\kappa}\, \frac1{\delta} \le \frac{C}\delta
\qquad \text{as}\ \delta\to 0\,,
\]
and note that the sequence
\[
m \mapsto a_m^2 r_0^{2m} = \frac{L(L+1)\, \ldots \, (L+m-1)}{m!}\, r_0^{2m}
\]
increases for $m\le [M]$ and decreases for $m\ge [M]+1$.  Thus, the maximal term of
this sequence does not exceed
\[
C M^{L-1} (1-\kappa\delta)^{2M} \le \frac{C}{\delta^{L-1}}\,,
\]
whence
\[
\lambda_m(\Sigma) \le \frac{C N}{\delta^{L-1}}\,, \qquad m=0, 1, \ldots , N-1\,.
\]
Also
\begin{multline*}
\det(\Sigma) = \prod_{m=0}^{N-1} \lambda_m(\Sigma)
\ge \prod_{m=0}^{N-1} Na_m^2 (1-\kappa\delta)^{2m}
\ge N^N (1-\kappa\delta)^{N(N-1)} \prod_{m=0}^{N-1}  c(m+1)^{L-1} \\
\ge (cN)^N (N!)^{L-1} (1-\kappa\delta)^{N(N-1)}
= \exp\bigl[ (1+o(1)) (LN\log N-\kappa\delta N^2) \bigr]
\end{multline*}
as $\delta\to 0$.

\subsection{Completing the proof}

Finally,
\begin{align*}
\log\bE\Bigl[ \prod_{j=1}^N |F(\omega^j\tau r_0)|^{-\theta}\Bigr] &\le -\log\det\Sigma + N\bigl( 1-\tfrac12\theta \Bigr) \log\Lambda
+ N \log\Gamma\bigl( 1-\tfrac12\theta \Bigr) \\
&\le -(1+o(1)) \bigl( LN\log N - \kappa\delta N^2 \bigr)  \\ \\
&\qquad+  N\bigl( 1-\tfrac12\theta \bigr)\, \bigl(\log N + (L-1)\log\tfrac1{\delta} + O(1)\bigr)
+ N \log\Gamma\bigl( 1-\tfrac12\theta \bigr)\,,
\end{align*}
and then,
\begin{align*}
\log \bP \Bigl[ &\sum_{j=1}^N \log|F(\omega^j \tau r_0 )| \le N \log |F(0)| + C \Bigr] \\
&\le N\theta \bigl( \tfrac12+a \bigr) \log\tfrac1{\delta} + O(1) + \log\bE\Bigl[ \prod_{j=1}^N |F(\omega^j\tau r_0)|^{-\theta}\Bigr] \\
&\le  N\theta \bigl( \tfrac12+a \bigr) \log\tfrac1{\delta} + O(1)
-(1+o(1)) \bigl( LN\log N - \kappa\delta N^2 \bigr)  \\ \\
&\qquad+  N\bigl( 1-\tfrac12\theta \bigr)\, \bigl(\log N + (L-1)\log\tfrac1{\delta} + O(1)\bigr)
+ N \log\Gamma\bigl( 1-\tfrac12\theta \bigr) \\ \\
& \stackrel{\rm def}{=} (1 + o(1)) N \cdot P_N.
\end{align*}
We set $\theta = 2 - a^2$, $\kappa = 1 + \delta$,
and continue to bound $P_N$ (using the choices \eqref{eq:choice_of_N} and \eqref{eq:choice_of_a}):
\begin{align*}
P_N & = (2-a^2)(\tfrac12 + a) \log \tfrac1\delta - L \log N +  (1+ \delta)\delta N + \tfrac{a^2}{2} (\log N + (L-1)\log \tfrac1\delta) + \log \Gamma(\tfrac{a^2}{2})\\ \\
& = \log \tfrac1\delta + O(\sqrt{\log \tfrac1\delta}) - L \log \tfrac1\delta + \tfrac{L-1}{2} \log \tfrac1\delta + O(\delta \log \tfrac1\delta) + O(\log \log \tfrac1\delta) + O(1) + O(\log \log \tfrac1\delta)\\ \\
& = -\tfrac{L-1}{2} \log \tfrac1\delta (1 + o(1)), \quad \delta\to 0.
\end{align*}

Thus,
\[
\log \bP \Bigl[ \sum_{j=1}^N \log|F(\omega^j \tau r_0 )| \le N \log |F(0)| + C \Bigr]
\le - (1+o(1))\, \frac{(L-1)^2}4\, \frac1{\delta}\, \Bigl( \log\frac1{\delta} \Bigr)^2\,.
\]
From Lemma~\ref{lemma9} we have
\begin{align*}
  \log\bP\bigl[ {\rm Hole}(r) \bigr] &\le C\log\frac1{\delta} + 2\log\frac1{\kappa - 1} - (1+o(1))\, \frac{(L-1)^2}4\, \frac1{\delta}\, \Bigl( \log\frac1{\delta} \Bigr)^2\\
  &=- (1+o(1))\, \frac{(L-1)^2}4\, \frac1{\delta}\, \Bigl( \log\frac1{\delta} \Bigr)^2\,,
\end{align*}
by our choice of $\kappa$, completing the proof of the upper bound in the case $L>1$. \hfill $\Box$

\section{Lower bound on the hole probability for $L>1$}

As before, let $\delta=1-r$ and assume that $r$ is sufficiently close to $1$. Introduce the parameter
\[
\tfrac12<\alpha<1\,,
\]
and put
\[
N = \Bigl[\, \frac{2L}{\delta}\, \log\frac1{\delta}\, \Bigr], \quad M = \frac1{\sqrt{\delta}}\, \Bigl( \log\frac1\delta\Bigr)^\alpha\,.
\]

We now introduce the events
\[
\cE_1 = \bigl\{ |\zeta_0|>M \bigr\}\,, \
\cE_2 = \Bigl\{ \max_{r\bT} \Bigl| \sum_{1\le n \le N } \zeta_n a_n z^n \Bigr| \le \frac{M}2 \Bigr\}\,, \
\cE_3 = \Bigl\{ \max_{r\bT} \Bigl| \sum_{n>N} \zeta_n a_n z^n \Bigr| \le \frac{M}2 \Bigr\}.
\]
Then
\[
\bP\bigl[ \text{Hole}(r)\bigr] \ge \bP\bigl[ \cE_1 \bigr] \cdot \bP\bigl[ \cE_2 \bigr] \cdot \bP\bigl[ \cE_3 \bigr]\,.
\]
We have
\[
\bP\bigl[ \cE_1 \bigr]  = \exp \bigl[ -M^2 \bigr]
= \exp \Bigl[ - \frac1\delta\, \Bigl( \log\frac1\delta\Bigr)^{2\alpha} \Bigr].
\]

In order to give a lower bound for $\bP\bigl[ \cE_2 \bigr]$, we rely on a comparison principle between Gaussian analytic functions which might be of use in other contexts. To introduce this principle let us say that a random analytic function $G$ has the $\mathrm{GAF}(b_n)$ distribution, for some sequence of complex numbers $(b_n)_{n\ge0}$, if $G$ has the same distribution as
\[z\mapsto\sum_{n\ge0} \zeta_n b_n z^n,\quad z\in\bC, \]
where, as usual, $(\zeta_n)$ is a sequence of independent standard complex Gaussian random variables.


\begin{lemma}\label{lemma_GAF_coupling}
Let $(b_n)$, $(c_n)$, $n\ge0$, be two sequences of complex numbers, such that $|c_n| \le |b_n|$ for all $n$, and put
\[
Q = \prod_{n\ge0} \left|\frac{c_n}{b_n}\right|,
\]
where we take the ratio $c_n / b_n$ to be $1$ if both $b_n$ and $c_n$ are zero. If $Q > 0$, then there exists a probability space supporting a random analytic function $G$ with the $\mathrm{GAF}(b_n)$ distribution, and an event $E$ satisfying
$\bP[E] = Q^2$, such that, conditioned on the event $E$, the function $G$ has the $\mathrm{GAF}(c_n)$ distribution.
\end{lemma}

The proof of Lemma \ref{lemma_GAF_coupling} uses the following simple property of Gaussian random variables.
\begin{lemma}\label{lemma_complex_Gaus_coupling}
Let $0 < \sigma \le 1$. There exists a probability space supporting a standard complex Gaussian random variable $\zeta$, and an event $E$ satisfying $\bP[E] = \sigma^2$, such that, conditioned on the event $E$, the random variable $\zeta$ has the complex Gaussian distribution with variance $\bE[|\zeta|^2\mid E] = \sigma^2$.
\end{lemma}
\noindent{\em Proof of Lemma~\ref{lemma_complex_Gaus_coupling}}:
We may assume that $\sigma < 1$. Write
$$ f(z) = \frac1\pi \exp(-|z|^2), \quad f_\sigma(z) = \frac1{\pi\sigma^2}\exp\big(-\tfrac1{\sigma^2} |z|^2 \big), \quad z \in \bC, $$
for the density of a standard complex Gaussian, and the density of a complex Gaussian with variance $\sigma^2$, respectively. Observe that, since $\sigma < 1$,
\begin{equation}\label{decomposition_Gaus_density}
f = \sigma^2 \cdot f_\sigma + (1-\sigma^2) g_\sigma
\end{equation}
for some non-negative function $g_\sigma$ with integral $1$. Now, suppose that our probability space supports a complex Gaussian $\zeta_\sigma$ with $\bE|\zeta_\sigma|^2 = \sigma^2$, a random variable $Y_\sigma$ with density function $g_\sigma$, and a Bernoulli random variable $I_\sigma$, satisfying $\bP[I_\sigma = 1] = 1 - \bP[I_\sigma~=~0] = \sigma^2$, which is independent of both $\zeta_\sigma$ and $Y_\sigma$. Then \eqref{decomposition_Gaus_density} implies that the random variable
\[\zeta = I_\sigma \cdot \zeta_\sigma + (1-I_\sigma) Y_\sigma \]
has the standard complex Gaussian distribution. In this probability space, after conditioning that $I_\sigma = 1$, the distribution of $\zeta$ is that of a complex Gaussian with variance $\bE[|\zeta|^2 | I_\sigma~=~1] = \sigma^2$, as required.
\hfill $\Box$

\medskip

\noindent{\em Proof of Lemma~\ref{lemma_GAF_coupling}}:
Let $\sigma_n = \left|c_n/b_n\right|$ (where again, the ratio is defined to be $1$ if $b_n$ and $c_n$ are both zero). Lemma \ref{lemma_complex_Gaus_coupling} yields, for each $n$, a probability space supporting a standard Gaussian random variable $\zeta_n$ and an event $E_n$ with $\bP[E_n] = \sigma_n^2$. Clearly we may assume that the sequence $\zeta_n$ is mutually independent and we take the probability space to be the product of these probability spaces, extend each $E_n$ to this space in the obvious way and define
$$G(z) = \sum_{n\ge0} \zeta_n b_n z^n.$$
The claim follows with the event $E=\bigcap_{n\ge0} E_n$.
\hfill $\Box$

\subsection{Estimating $ \bP\bigl[ \cE_2 \bigr] $}
Put
\[ G(z) = \sum_{1\le n \le N } \zeta_n a_n z^n, \]
and let $(q_n)_{n=1}^N$ be a sequence of numbers in $[0,1]$ to be specified below. According to Lemma \ref{lemma_GAF_coupling}, there is an event $E$ with probability $Q^2 = \prod_{n=1}^N q_n^2$, such that on the event $E$, the function $G$ has the same distribution as
\[ G_Q(z) \stackrel{\rm def}= \sum_{1\le n \le N } \zeta_n q_n a_n z^n. \]
Notice also that
\[
\sigma_Q(\rho)^2 \stackrel{\rm def}{=} \bE \Bigl[ \Bigl| \sum_{n=1}^N \zeta_n q_n a_n (\rho e^{{\rm i}\theta})^n \Bigr|^2 \Bigr] = \sum_{n=1}^N q_n^2 a_n^2 \rho^{2n}\,.
\]
If we now set
\[
r_2 = r + \delta^2, \quad
\lambda = \Big(\log \frac1\delta\Big)^\alpha,
\]
then applying first Lemma \ref{lemma_GAF_coupling}, and then Lemma \ref{lemma2} to the function $G_Q$, with $\la$ as above, we obtain that for $\delta$ sufficiently small
\begin{align*}
\bP \Bigl[ \max_{r \bar\bD} \Bigl| \sum_{n=1}^N \zeta_n a_n z^n \Bigr| \le \la \, \sigma_Q(r_2) \Bigr]
& \ge
Q^2 \cdot \bP \Bigl[ \max_{r \bar\bD} \Bigl| \sum_{n=1}^N \zeta_n a_n z^n \Bigr| \le \la \, \sigma_Q(r_2) \big\rvert E \Bigr] \\
& =
Q^2 \cdot \bP \Bigl[ \max_{r \bar\bD} \Bigl| \sum_{n=1}^N \zeta_n q_n a_n z^n \Bigr| \le \la \, \sigma_Q(r_2) \Bigr]\\
&\ge Q^2 \Bigl(1 - C \delta^{-1} \exp\bigl(-c (\log \tfrac{1}\delta )^{2\alpha} \bigr) \Bigl)\\
& \overset{\alpha > \tfrac12}{\ge} \tfrac12 Q^2.
\end{align*}

%

For this estimate to be useful in our context we have to ensure, by choosing the sequence $q_n$ appropriately, that
\begin{equation}\label{eq:cond_for_sigma_q}
\sigma_Q(r_2) \le \frac{M}{2 \la} = \frac1{2 \sqrt{\delta}}.
\end{equation}
First, it is straightforward to verify that
\[ \exp(-\delta) \le r_2 \le \exp(-\delta +\delta^2), \quad \forall r \in (0,1). \]
Since $c n^{L-1} \le a_n^2 \le C n^{L-1}$, this implies that for $n \in \{1, \dots, N \}$
\[
a_n^2 r_2^{2 n} = \exp\big((L-1)\log n - 2 n \delta + O(1) \big).
\]
Putting
\[
N_1 = \Bigl[\, \frac{L - 1}{2 \delta}\, \log\frac1{\delta}\, \Bigr]
\]
and noting that the function $x \mapsto (L-1)\log x - 2 x \delta$ attains its maximum at $x = \frac{L-1}{2 \delta}$ we see that $a_{n}^2 r_2^{2 n} \ge c$ for $n \le N_1$, while for $n \in \{N_1, \dots, N\}$ we have $a_{n}^2 r_2^{2 n} \le C \left(\log \frac1\delta \right)^{L-1}$.

We therefore choose
\[
q^2_n = \begin{cases}
\alpha_1 \big( a^2_n r_2^{2n} \log \tfrac1\delta \big)^{-1}\,, &  n\in\{1, \dots, N_1\},\\ \\
\alpha_1 \big( \log \tfrac1\delta \big)^{-L}\,, &  n\in\{N_1 + 1, \dots, N\},
\end{cases}
\]
where we choose the constant $\alpha_1>0$ sufficiently small to ensure that $q_n \le 1$ for all $n \in \{1, \dots N\}$. With this choice we have
\begin{align*}
\sigma^2_Q(r_2) = \sum_{n=1}^N q^2_n a^2_n r_2^{2n} \le \alpha_1 \Big(\log \frac1\delta\Big)^{-1} N_1 + \alpha_1 \cdot C \Big(\log \frac1\delta\Big)^{-1} (N - N_1) \le C \cdot \frac{\alpha_1}{\delta}
\end{align*}
and further choosing $\alpha_1$ small if necessary, Condition \eqref{eq:cond_for_sigma_q} is satisfied.

It remains to estimate the probability of the event $E$. Notice that
\begin{align*}
\bP[E] = Q^2 & = \prod_{n=1}^N q^2_n = \alpha_1^N \big( \log \tfrac1\delta \big)^{-L(N - N_1) - N_1} \prod_{n=1}^{N_1} \big( a^2_n r_2^{2n}\big)^{-1}\\
& \ge \exp\big( -C N \log \log \tfrac1\delta \big) \cdot \prod_{n=1}^{N_1} \big( C n^{L-1} \exp[-2n(\delta - \delta^2)] \big)^{-1} \\
& \ge \exp\big( -C N \log \log \tfrac1\delta \big) \cdot \frac{e^{(\delta - \delta^2 )N_1(N_1+1)}}{(N_1!)^{L-1}}.
\end{align*}
Recalling that $N \le C N_1$, and using that $N_1 = \bigl[\, \tfrac{L - 1}{2 \delta}\, \log\tfrac1{\delta}\, \bigr]$, we obtain
\begin{align*}
\bP\bigl[ \cE_2 \bigr] & \ge \tfrac12 Q^2 \ge \exp\big(- N_1 \bigl[\, (L-1) \log N_1 - \delta N_1 + C \log\log \tfrac1\delta \,\bigr]\big)\\ \\
& = \exp\big(- \tfrac14 [(L-1)^2 + o(1)]\, \tfrac1\delta \, \log^2\tfrac1\delta \big).
\end{align*}

\subsection{Estimating $ \bP\bigl[ \cE_3 \bigr] $}

Put
\[
H(z) = \sum_{n>N} \zeta_n a_n z^n\,.
\]
Then
\[
\sigma_H(\rho)^2 = \bE \Bigl[ \Bigl| \sum_{n>N} \zeta_n a_n (\rho e^{{\rm i}\theta})^n \Bigr|^2 \Bigr] = \sum_{n>N} a_n^2 \rho^{2n}\,.
\]
The choice of the parameter $N$ guarantees that the `tail' $H$ will be small.

\begin{lemma}\label{lemma_tail_bnd_Lg1}
Put $r_1=\tfrac12 (1+r)$.
There exists a constant $C > 0$ such that,
for $\delta$ sufficiently small,
\[
\sigma_H(r_1)^2 = \sum_{n>N} a_n^2 r_1^{2n} \le C.
\]
\end{lemma}
\noindent{\em Proof of Lemma~\ref{lemma_tail_bnd_Lg1}}:
Since the function $r \mapsto 2 \log(1+r) - r$, $0\le r \le 1$ attains its maximum at $r=1$ we have
$r_1^2 \le \exp(-\delta)$ for $\delta = 1 - r \in [0,1]$. Recalling that $N = \bigl[\, \tfrac{2L}{\delta}\, \log\tfrac1{\delta}\, \bigr]$, we observe that for $n > N$ we have
\[
\frac{n}{\log n} \ge \frac{2(L-1)}{\delta},
\]
and therefore $n^{L-1} \le \exp(\tfrac12 \delta n)$.
Then, using that $a_n^2 \le C n^{L-1}$, we have
\[
\sigma_H(r_1)^2 \le C \sum_{n>N} n^{L-1} \exp(-\delta n)
 \le C \sum_{n>N} \exp(-\tfrac12 \delta n) \le C \delta^{-1} \exp(-\tfrac12 \delta N) \le C \delta^{L - 1}.
\]
Since $L>1$ this is a stronger result than we claimed.
\hfill $\Box$

\medskip
By Lemma~\ref{lemma_tail_bnd_Lg1},
\[
\bP\bigl[ \cE_3^c \bigr]= \bP \bigl[ \max_{r\bar\bD} |H| > \tfrac12 M \bigr]
\le \bP \bigl[ \max_{r\bar\bD} |H| > c M \cdot \sigma_{H}\bigl( \tfrac12 (1+r) \bigr) \bigr]\,,
\]
if $c>0$ is sufficiently small. By Lemma~\ref{lemma2}, the right-hand side is at most
$C\delta^{-1}\, \exp\bigl[ -c M^2 \bigr] \to 0 $ as
$ \delta\to 0 $.
Thus, $\bP\bigl[ \cE_3 \bigr] \ge \tfrac12 $ for $\delta$ sufficiently small.

\subsection{Putting the estimates together}
Finally,
\begin{align*}
\bP \bigl[ \text{Hole}(r) \bigr] &\ge \bP \bigl[ \cE_1 \bigr] \cdot
\bP \bigl[ \cE_2 \bigr] \cdot \bP \bigl[ \cE_3 \bigr] \\ \\
&\ge \frac12\, \exp\biggl[ -\frac1\delta\, \bigl( \log\frac1\delta \bigr)^{2\alpha} -
\frac{(L-1)^2+o(1)}4\, \frac1\delta\, \log^2\frac1{\delta} \biggr] \\ \\
&\stackrel{\alpha<1}\ge \exp\biggl[
- \frac{(L-1)^2+o(1)}4\, \frac1\delta\, \log^2\frac1{\delta} \biggr]\,,
\end{align*}
completing the proof of the lower bound in the case $L>1$, and hence of Theorem~\ref{thm:main}.
\hfill $\Box$

\section{The hole probability for non-invariant GAFs with regularly distributed coefficients}\label{Section:Last}

The proofs we gave do not use the hyperbolic invariance of the zero distribution of $F=F_L$.
In the case $0<L<1$ we could assume that the sequence of coefficients $(a_n)$ in \eqref{eq1} does not increase, that $a_0=1$ and that
\begin{equation}\label{eq:coefficients}
a_n \simeq n^{\frac12 (L-1)}\,,\quad n\ge 1.
\end{equation}
Then, setting as before
\begin{equation*}
  \sigma_F(r)^2 = \bE\bigl[ |F(re^{{\rm i}\theta})|^2\bigr] = \sum_{n\ge 0} a_n^2 r^{2n}
\end{equation*}
the same proof yields the bounds
\begin{align*}
\frac{1 - L + o(1)}2\, \sigma_F(r)^2 \log \frac1{1-r} &\le
- \log \bP\bigl[ \text{Hole}(r) \bigr] \\
&\le (1 - L + o(1))\, \sigma_F(r)^2 \log \frac1{1-r}\,, \qquad r\to 1\,.
\end{align*}
In the case $L>1$, assuming that $a_0=1$ and \eqref{eq:coefficients}, we also get the same answer as in Theorem~\ref{thm:main}:
\[
- \log \bP\bigl[ \text{Hole}(r) \bigr]
= \frac{(L-1)^2 +o(1)}4\, \frac1{1-r}\, \log^2\frac1{1-r}\,,
\qquad r\to 1\,.
\]

In the case $L=1$ (that is, $a_n=1$, $n\ge 0$)
the result of Peres and Vir\'ag relies on their proof that the zero set of $F_1$ is a determinantal point process.
This ceases to hold under a slight perturbation of the coefficients $a_n$,
while our techniques still work, though yielding less precise bounds:

\begin{theorem}\label{Thm:L=1}
Suppose that $F$ is a Gaussian Taylor series of the form \eqref{eq1} with $a_n \simeq 1$ for $n\ge 0$. Then
\[
- \log \bP\bigl[ {\rm Hole}(r) \bigr] \simeq \frac{1}{1-r}\,, \qquad \frac12 \le r < 1\,.
\]
\end{theorem}

For the reader's convenience, we supply the proof of Theorem~\ref{Thm:L=1}, which is based on arguments similar
to those we have used above.
As before, we put $\delta=1-r$, $\sigma_F=\sigma_F(r)$,
and note that under the assumptions of Theorem~\ref{Thm:L=1},
$ \sigma_F(r)^2 \simeq \delta^{-1}$. Also put
\begin{equation*}
N=\left\lceil\frac{1}{1-r}\right\rceil
\end{equation*}
and $\omega = e(1/N)$.

\subsection{Upper bound on the hole probability in
Theorem~\ref{Thm:L=1}}\label{sec:upper_bound_L_1}

\subsubsection{Beginning the proof}
The starting point is the same as in the proofs of the upper bound in the cases $L\neq 1$.
Put $r_0 = 1-2\delta$. Then, by Lemma~\ref{lemma9} (with $\kappa=2$ and $k=N$),
\[
\bP\bigl[ \text{Hole}(r)  \bigr]
\le\, \delta^{-c}\,
\sup_{\tau\in\bT}\, \bP\biggl[ \sum_{j=1}^N \log|F(\omega^j \tau r_0 )| \le N\log |F(0)| + C \biggr]
+ \langle {\rm negligible\ terms} \rangle\,.
\]
For fixed $\tau\in\bT$ and for a small positive parameter $b$, we have
\begin{align*}
\bP\biggl[ &\sum_{j=1}^N \log|F(\omega^j \tau r_0 )| \le N\log |F(0)| + C \biggr] \\
&\le \bP\biggl[ \sum_{j=1}^N \log|F(\omega^j \tau r_0 )| \le N\log\bigl( b\sigma_F \bigr) + C \biggr]
+ \bP\bigl[ |F(0)|>b\sigma_F \bigr] \\
&\le \bP\biggl[ \sum_{j=1}^N \log |w_j| \le N\log\bigl( C b \bigr)\biggr]
+ \bP\bigl[ |F(0)|>b\sigma_F \bigr]\,,
\end{align*}
where $w_j = F(\omega^j \tau r_0)/\sqrt{N}$ are complex Gaussian random variables.
Next,
\[
\bP\Bigl[ \sum_{j=1}^N \log |w_j| \le N\log\bigl( C b \bigr)\Bigr] = \bP \Bigl[ \prod_{j=1}^N  |w_j|^{-1} \ge (Cb)^{-N} \Bigr]
\le (Cb)^N\, \bE \Bigl[ \prod_{j=1}^N  |w_j|^{-1} \Bigr]\,.
\]
Thus, up to negligible terms, the hole probability is bounded from above by
\[
(Cb)^N\, \bE \Bigl[ \prod_{j=1}^N  |w_j|^{-1} \Bigr] + e^{-b^2 \sigma_F^2}\,.
\]
What remains is to show that the expectation of the product of $|w_j|^{-1}$ grows at most exponentially
with $N$. Then, choosing the constant $b$ so small that the prefactor  $ (Cb)^N $ overcomes this growth of
the expectation, we will get the result.

\subsubsection{Estimating $ \bE \Bigl[ \prod_{j=1}^N  |w_j|^{-1} \Bigr] $}

One can use Lemma~\ref{lemma18} in order to bound the expectation above, below we give an alternative argument. Put $z_j = \omega^j \tau r_0$, $1\le j \le N$, and consider the
covariance matrix
\[
\Gamma_{ij} = \bE \bigl[ w_i \bar w_j\bigr] = N^{-1} \bE \bigl[ F(z_i) \bar F(z_j)\bigr]\,,
\qquad 1 \le i, j \le N\,.
\]
For each non-empty subset $I\subset \{1, 2, \ldots , N\}$, we put
$ \Gamma_I = \bigl( \Gamma_{ij} \bigr)_{i, j\in I}$.

\begin{lemma}\label{lemma:minors}
For each $I\subset \{1, 2, \ldots , N\}$, we have
$ \det \Gamma_I \ge c^{|I|} $.
\end{lemma}

\noindent{\em Proof of Lemma~\ref{lemma:minors}}:
By Lemma~\ref{lemma14}, the eigenvalues of the matrix $\Gamma$ are
\[
\la_m = \sum_{n\equiv m\, (N)}\, a_n^2 r_0^{2n} \ge c \sum_{n\equiv m\, (N)} (1-2\delta)^{2n}
\ge c \sum_{k\ge 0} e^{-CNk\delta} \ge c \sum_{k\ge 0} e^{-Ck} \ge c >0\,,
\]
that is, the minimal eigenvalue of $\Gamma$ is separated from zero.
It remains to recall that the $N-1$ eigenvalues of any minor of order $N-1$ of an Hermitian matrix  of order $N$ interlace with the $N$ eigenvalues of the original matrix.
Applying this principle several times, we conclude that the minimal eigenvalues of the matrix
$\Gamma_I$ cannot be less than the minimal eigenvalue of the full matrix $\Gamma$. \hfill $\Box$

\medskip
Now, we write
\[
 \bE \biggl[ \prod_{j=1}^N  |w_j|^{-1} \biggr]
 \le \sum_{I\subset\{1, 2, \ldots , N\}}\,  \bE \biggl[ \prod_{i \in I}  |w_i|^{-1} \done_{\{|w_i|\le 1\}} \biggr]\,.
\]
By Lemma~\ref{lemma:minors}, the expectations on the right-hand side do not exceed
\[
\frac1{\pi^{|I|} \det \Gamma_I}\,
\biggl( \int_{|w|\le 1} \frac1{|w|}\, e^{-\la_I^{-1}|w|^2 }\, {\rm d}m(w) \biggr)^{|I|}
\le C^{|I|}
\]
(here, $\la_I$ is the maximal eigenvalue of $\Gamma_I$).
Since the number of subsets $I$ of the set $ \bigl\{ 1, 2, \ldots , N \bigr\} $ is $2^N$, we finally
get
\[
 \bE \biggl[ \prod_{j=1}^N  |w_j|^{-1} \biggr] \le (2C)^N\,.
\]
This completes the proof of the upper bound on the hole probability in Theorem~\ref{Thm:L=1}. \hfill $\Box$

\subsection{Lower bound on the hole probability in Theorem~\ref{Thm:L=1}}
We write
\[
F(z) = F(0) + G(z) =  F(0) + \sum_{k \ge 0} z^{k N} S_k(z)\,,
\]
where the $S_k (z)$ are independent random Gaussian polynomials,
\[
S_k(z) = \sum_{j=1}^N \zeta_{j+kN} a_{j+kN} z^j\,.
\]
We have
\[
\max_{r\bT} |G| \le \sum_{k\ge 0} r^{kN} \max_{r\bT} |S_k| \, \stackrel{N\log r < -1}\le \,
\sum_{k\ge 0} e^{-k}  \max_{r\bT} |S_k|\,.
\]
We fix $1<A<e$ and consider the independent events
$ \cE_k = \bigl\{ \max_{r\bT} |S_k| < A^k \sqrt{N} \bigr\} $. If these events occur together,
then
\[
\max_{r\bT} |G| < \sqrt{N}\, \sum_{k\ge 0} \bigl( Ae^{-1} \bigr)^k < B\sqrt{N}
\]
with some positive numerical constant $B$. Then, $F(z) \ne 0$ on $r\bar\bD$, provided that
$|F(0)| > B\sqrt{N}$, and that all the events $\cE_k$ occur together, that is,
\[
\bP\bigl[ \text{Hole} (r) \bigr] \ge e^{-B^2 N}\, \prod_{k\ge 0} \bP \bigl[ \cE_k \bigr]\,.
\]
It remains to estimate from below the probability of each event $\cE_k$ and to multiply the estimates.

\subsubsection{Estimating $  \bP \bigl[ \cE_k \bigr] $}

Take $4N$ points $z_j = r e\big(\frac{j}{4N}\big)$, $1\le j \le 4N$. By the Bernstein inequality,
\[
\max_{\theta\in [0, 2\pi]} \Bigl| \frac{\partial S_k}{\partial\theta} \bigl( re^{{\rm i}\theta} \bigr)\Bigr|
\le N\, \max_{r\bT} |S_k|\,.
\]
Therefore,
\[
\max_{r\bT} |S_k| \le \max_{1\le j \le 4N}\, |S_k(z_j)| + \frac{\pi}{4N}\cdot N \max_{r\bT} |S_k|\,,
\]
whence,
\[
\max_{r\bT} |S_k| \le C \max_{1\le j \le 4N}\, |S_k(z_j)|\,,
\]
and
\[
\bP \bigl[ \cE_k \bigr] \ge
\bP\Bigl[\, \max_{1\le j \le 4N}\, |S_k(z_j)| < C A^k \sqrt{N} \, \Bigr]\,.
\]
Then, applying as in Section~\ref{section:LB,L<1}, Harg\'e's version of the
Gaussian correlation inequality, we get
\[
\bP \bigl[ \cE_k \bigr] \ge
\prod_{j=1}^{4N} \bP\bigl[ |S_k(z_j)| < C A^k \sqrt{N} \bigr]\,.
\]
(In fact, passing to real and imaginary parts we could use here the simpler Khatri-Sidak version of the Gaussian correlation
inequality~\cite[Section~2.4]{LS}.)
Each value $S_k(z_j)$ is a complex Gaussian random variable with variance
\[
\sigma_{S_k}^2 = \sum_{j=1}^N a_{j+kN}^2 r^{2j} \simeq\, \sum_{j=1}^N r^{2j} \simeq \frac1{\delta}
\simeq N\,.
\]
Therefore,
\[
 \bP\bigl[ \, |S_k(z_j)| < C A^k \sqrt{N}\, \bigr] =
 1 -  \bP\bigl[\, |S_k(z_j)| \ge C A^k \sqrt{N} \, \bigr]
 \ge 1 - e^{-cA^{2k}}\,,
\]
whence,
\[
\bP \bigl[ \cE_k \bigr]
\ge \Bigl( 1 - e^{-cA^{2k}} \Bigr)^{4N}\,,
\]
and then,
\[
\prod_{k\ge 0} \bP \bigl[ \cE_k \bigr]
\ge \prod_{k\ge 0} \Bigl( 1 - e^{-cA^{2k}} \Bigr)^{4N} \ge e^{-cN}\,,
\]
completing the proof. \hfill $\Box$

\end{document}